\documentclass[11pt]{article}
\usepackage{amssymb,amsmath,amsthm,amsfonts}
\usepackage[dvips]{graphicx, color}
\usepackage{subfigure}
\usepackage{latexsym}

\theoremstyle{plain}
\newtheorem{theorem}{Theorem}[section]
\newtheorem{lemma}{Lemma}[section]
\newtheorem{proposition}{Proposition}[section]

\numberwithin{equation}{section}
\theoremstyle{definition}

\newtheorem{remark}{Remark}[section]

\begin{document}

\title{\textbf{Local Exact Controllability of a Parabolic System of
Chemotaxis}\thanks{%
This work was supported by Program for Innovative Research Team in UIBE, the
National Natural Science Foundation of China, the National Basic Research
Program of China (2011CB808002), the National Research Foundation of South
Africa, the National Science Foundation of China (11201358), and the
Fundamental Research Funds for the Central Universities.}}
\date{December 23, 2012}
\author{Bao-Zhu Guo$^{a,c}$ ~and Liang Zhang$^{b,c}$\thanks{%
The corresponding author. Email: changleang@yahoo.com.cn} \\
$^a$\textit{Academy of Mathematics and Systems Science, Academia Sinica}\\
\textit{Beijing 100190, China} \\
$^b$\textit{Department of Mathematics,}\\
\textit{Wuhan University of Technology, Wuhan 430070, China}\\
$^c$\textit{School of Computational and Applied Mathematics}\\
\textit{University of the Witwatersrand, Wits 2050, Johannesburg, South
Africa}}
\maketitle

\begin{abstract}
This paper studies  the controllability problem of a parabolic
system of chemotaxis. The local exact controllability to
trajectories of the system imposed  one control force only is
obtained by applying Kakutani's fixed point theorem combined with
the null controllability of the associated linearized parabolic
system. The control function is shown to be in $L^\infty(Q)$, which
is estimated by using the methods of maximal regularity and
$L^p$-$L^q$ estimates of parabolic equations.

\vspace{0.3cm}

\textbf{Keywords:}~local exact controllability, chemotaxis system, Carleman
inequality, Kakutani's fixed point theorem.

\vspace{0.3cm}

\textbf{AMS subject classifications:}~ 93B05, 93C20, 35B37.
\end{abstract}

\section{Introduction and main results}

Let $\Omega \subset {\mathbb{R}}^{N}(N\geq 1)$ be a bounded domain with
sufficient smooth boundary $\partial \Omega $. Let $\omega $ be a nonempty
open subset of $\Omega $, and $T>0$. We denote $Q=\Omega \times (0,T)$, $%
\Sigma =\partial \Omega \times (0,T)$ and $Q_{\omega }=\omega \times (0,T)$.
Throughout this paper, we use $W^{s,q}(\Omega )$, $W_{q}^{2,1}(Q)$ and $%
C^{\alpha }(\overline{\Omega })$ ($s,\alpha \geq 0$, $1\leq q\leq
\infty $) for the usual Sobolev spaces (e.g., \cite{lady}), and set
$H^{m}(\Omega )=W^{m,2}(\Omega )$ for $m\in {\mathbb{N}}$.
$L^{p}(\Omega )$ and $L^{p}(Q)$ ($1\leq p\leq \infty $) are the
usual Lebesgue function spaces with the norm $\vert \cdot \vert
_{p}$ and $\Vert \cdot \Vert _{p}$, respectively. Moreover, let
\begin{eqnarray*}
V^{1}(Q) &=&\left\{ y|y\in L^{2}(0,T;H^{1}(\Omega )),\partial _{t}y\in
L^{2}(0,T;H^{1}(\Omega )^{\ast })\right\} , \\
V^{2}(Q) &=&\left\{ y|y\in L^{2}(0,T;H^{2}(\Omega )),\partial _{t}y\in
L^{2}(Q)\right\} ,
\end{eqnarray*}%
be equipped with their graph norms, where $H^{1}(\Omega )^{\ast }$
denotes the dual space of $H^{1}(\Omega )$. The duality between
$H^{1}(\Omega )^{\ast }$ and $H^{1}(\Omega )$ is denoted by $\langle
\cdot ,\cdot \rangle .$

\bigskip

In this paper, we are concerned with the following controlled parabolic
system with state functions $u\equiv u(x,t)$ and $v\equiv v(x,t):$
\begin{equation}
\begin{cases}
\partial _{t}u=\nabla \cdot \left( \nabla u-\chi u\nabla v\right) +{\mathbf{1%
}}_{\omega }f & \mathrm{in}\ Q, \\
\partial _{t}v=\Delta v-\gamma v+\delta u & \mathrm{in}\ Q, \\
{\partial }_{\nu }u=0,\partial _{\nu }v=0 & \mathrm{on}\ \Sigma , \\
u(x,0)=u_{0}(x)\ \ v(x,0)=v_{0}(x) & x\in \Omega ,%
\end{cases}
\label{e}
\end{equation}%
where $\partial _{t}=\partial /\partial t$, and $\partial _{\nu }=\partial
/\partial \nu $ stands for the derivative with respect to the outer normal $%
\nu $ of $\partial \Omega $, $\mathbf{1}_{\omega }$ represents the
characteristic function of $\omega $, $f\equiv f(x,t)$ is the control
function so that $\mathbf{1}_{\omega }f$ is the control force acting from
the outside on a portion of the domain $\Omega $, $u_{0}$ and $v_{0}$ are
the initial values, and $\chi ,\gamma $ and $\delta $ are given positive
constants.

A pair of functions $(u,v)$ with
\begin{equation*}
u\in V^{1}(Q)\cap L^{\infty }(Q),v\in V^{2}(Q)\cap L^{\infty }(Q)
\end{equation*}%
is called a \textit{weak solution of \eqref{e}} if for all $\varphi \in
L^{2}(0,T;H^{1}(\Omega ))$, the following identities hold:%
\begin{eqnarray*}
\int_{0}^{T}\left\langle \partial _{t}u,\varphi \right\rangle dt+\iint_{Q}
\left[ (\nabla u-\chi u\nabla v)\cdot \nabla \varphi +{\mathbf{1}}_{\omega
}f\varphi \right] dxdt &=&0, \\
\iint_{Q}\varphi \partial _{t}vdxdt+\iint_{Q}\left[ \nabla v\cdot \nabla
\varphi +(-\gamma v+\delta u)\varphi \right] dxdt &=&0.
\end{eqnarray*}

We write the free system of \eqref{e}, that is, in the absence of $f$, as
follows:
\begin{equation}
\begin{cases}
\partial _{t}\overline{u}=\nabla \cdot \left( \nabla \overline{u}-\chi
\overline{u}\nabla \overline{v}\right) & \mathrm{in}\ Q, \\
\partial _{t}\overline{v}=\Delta \overline{v}-\gamma \overline{v}+\delta
\overline{u} & \mathrm{in}\ Q, \\
{\partial }_{\nu }\overline{u}=0,\partial _{\nu }\overline{v}=0 & \mathrm{on}%
\ \Sigma , \\
\overline{u}(x,0)=\overline{u}_{0}(x)\ \ \overline{v}(x,0)=\overline{v}%
_{0}(x) & x\in \Omega .%
\end{cases}
\label{ks}
\end{equation}%
The system (\ref{ks}) is a prototype chemotaxis system so called
Keller-Segel model which describes the aggregation process of slime
mold resulting from chemotactic attraction. In \eqref{ks},
$\overline{u}$ represents the density of the cellular slime mold,
$\overline{v}$ is the density of the chemical substance (see
\cite{keller}). In the last decade, there is a large number of works
devoted to the  mathematical analysis of the Keller-Segel system.
Several topics on the Keller-Segel model for chemotaxis such as
aggregation, blow-up of solutions, and chemotactic collapse, etc.,
have been concerned and some significant results
have been achieved from different discipline perspectives. In Horstmann \cite%
{horstmann1} and Hillen and Painter \cite{hillen}, it provides a
detailed introduction into the mathematics of the Keller-Segel model
for chemotaxis with abundant references therein. Here we would
mention a few facts about the local and global existence of
solutions for the Keller-Segel model. Generally speaking, the
blow-up of solutions of Keller-Segel system in finite or infinite
time depends strongly on the space dimension. In 1-d case, a finite
time blow-up never occur, and the global solution exists and
converges to the stationary solution as times goes to infinity (see \cite%
{osaki}). But the blow-up may occur in finite or infinite time in $n$%
-dimensional case for $n\geq 3$ (see \cite{boy,horstmann}). For the 2-d
case, several thresholds have been found. When the mass of the initial data
is below some threshold value, the solution exists globally in time and its $%
L^{\infty }$-norm is uniformly bounded for all time. While the mass of the
initial data is larger than some threshold value, the solution will blow up
either in finite or in infinite time (see \cite{biler,gajewski,yagi2}).

\bigskip Due to blow-up feature of solutions of the Keller-Segel model, it
is interesting to consider some controllability problems. Let $(\overline{u},%
\overline{v})$ be a trajectory, i.e., a solution of \eqref{ks} corresponding
to some initial value $(\overline{u}_{0},\overline{v}_{0})$. We say that the
system \eqref{ks} is \textit{locally exactly controllable to the trajectory }%
$(\overline{u},\overline{v})$\textit{\ at time }$T$, if there exists a
neighborhood $\mathcal{O}$ of $(\overline{u}_{0},\overline{v}_{0})$ such
that for any initial data $(u_{0},v_{0})\in \mathcal{O}$, the solution $%
(u,v) $ of \eqref{e} driven by some control function $f$ satisfies
\begin{equation*}
u(x,T)=\overline{u}(x,T),v(x,T)=\overline{v}(x,T),\hbox{ for}\ x\in \Omega \ %
\hbox{ a.e.},
\end{equation*}%
where the neighborhood $\mathcal{O}$ and the control function space will be
specified later.

\bigskip

In this paper, we suppose that $\overline{u},\overline{v}$ verify the
following regularity properties:
\begin{equation}
\overline{u},\overline{v}\in L^{\infty }(Q),\nabla \overline{v}\in L^{\infty
}(Q)^{N}.  \label{regularity}
\end{equation}

\begin{remark}
\textrm{The solution $(\overline{u},\overline{v})$ of the system
\eqref{ks} exists at least locally in time interval $[0,T_{1}]$ for
$T_{1}<T_{max}$ with sufficiently small initial data
$(\overline{u}_{0},\overline{v}_{0})$, where $T_{max}$ is the
maximal existence time (see \cite{horstmann1} and  reference
therein). If the system is locally exactly controllable, then we can
drive
the state of the system by some control force to a given trajectory at time $%
T\leq T_{1}$ before the time $T_{max}$ to avoid blow-up. It is also worth
indicating that the reason we consider the local exact controllability
instead of exact controllability is that the solution may blow up when the
mass of initial value is larger than some threshold value. }
\end{remark}

\begin{remark}
\textrm{When $(\overline u,\overline v)=(0,0)$, the local exact
controllability is reduced to the local null controllability. If the system %
\eqref{e} is locally null controllable at time $T$ with some control, then
we can switch off the control after time $T$ and the system will keep into
zero afterwards. }
\end{remark}

\bigskip

This paper is devoted to the local exact controllability of the coupled
parabolic system \eqref{e} via one control. The controllability of parabolic
systems of coupled equations attracts intensive attention in the last few
years. In Barbu \cite{barbu}, it studies the local exact controllability to
steady states with controls acting on each equation of the system via the
same interior domain. This could be done by taking it as a direct
consequence of the controllability of the scalar parabolic equations. It is
much more interesting and applicable to consider the controllability of a
parabolic system with one control force imposed on one equation of the
system. Ammar Kdjodia \textit{et al.~}\cite{ammar} is the first work of this
kind. They show that the phase-field system is locally exactly controllable
to the trajectory by one control force. The series of works of Ammar Kdjodia
\textit{et al.~}(\cite{ammar2,ammar3,ammar4}), and the works of Gonz\'{a}%
lez-Burgos \textit{et al.} (\cite{gonzalez,gonzalez1,teresa,fernandez2}),
have extended such problem to more general cases. The survey paper \cite%
{ammar0} gives a comprehensive introduction to this topic. For more works of
the controllability of parabolic equations, we also refer to \cite%
{fernandez1,fernandez3,fernandez4,fursikov} and \cite{coron}.

\bigskip

However, as to our best knowledge, very few results are available to the
control problems of the system \eqref{e}. In Ryu and Yagi \cite{yagi1}, it
considers an optimal control problem of the system \eqref{e} with the
control to be distributed on the second equation of \eqref{e}. The present
paper can be considered as a first work on the controllability of the system %
\eqref{e}. There are some other kinds of interesting control
problems for the system \eqref{e}. In the system \eqref{e}, the
chemotactic term $-\chi \nabla \cdot (u\nabla v)$ causes much more
mathematical difficulties than the coupled parabolic systems
aforementioned. The techniques presented in this paper would be
useful for other forms of chemotaxis system such as the
parabolic-elliptic chemotaxis system, and even for other coupled
systems like drift-diffusion equations from the semiconductor
device.

\bigskip

The idea of obtaining the controllability of \eqref{e} is somehow
classical: We first establish the null controllability of the
linearized system and then apply the fixed point theorem. Now we
consider the null controllability of the linearized system of
\eqref{e}, which is written as follows:
\begin{equation}
\begin{cases}
\partial _{t}y=\Delta y-\nabla \cdot \left( By\right) -\nabla \cdot \left(
a\nabla z\right) +\mathbf{1}_{\omega }f & \mathrm{in}\ Q, \\
\partial _{t}z=\Delta z-\gamma z+\delta y & \mathrm{in}\ Q, \\
{\partial }_{\nu }y=0,\partial _{\nu }z=0 & \mathrm{on}\ \Sigma , \\
y(x,0)=y_{0}(x)\ \ z(x,0)=z_{0}(x) & x\in \Omega ,%
\end{cases}
\label{lineare}
\end{equation}%
where $a\in L^{\infty }(Q)$, $B\in L^{\infty }(Q)^{N}$ with $B\cdot \nu =0$
on $\Sigma $, $f\in L^{2}(Q)$ is the control force, and $y_{0},z_{0}\in
L^{2}(\Omega )$ are given initial data. To study the null controllability of %
\eqref{lineare}, we are led to consider the observability of the adjoint
system of \eqref{lineare}:%
\begin{equation}
\begin{cases}
-\partial _{t}\phi =\Delta \phi +B\nabla \phi +\delta \theta & \mathrm{in}\
Q, \\
-\partial _{t}\theta =\Delta \theta -\gamma \theta -\nabla \cdot \left(
a\nabla \phi \right) & \mathrm{in}\ Q, \\
{\partial }_{\nu }\phi =0,\partial _{\nu }\theta =0 & \mathrm{on}\ \Sigma ,
\\
\phi (x,T)=\phi ^{T}(x),\theta (x,T)=\theta ^{T}(x) & x\in \Omega ,%
\end{cases}
\label{adjoint0}
\end{equation}%
where $\phi ^{T},\theta ^{T}\in L^{2}(\Omega )$. It is well-known that the
null controllability of \eqref{lineare} is equivalent to the observability
inequality for system (\ref{adjoint0}):%
\begin{equation*}
\left\vert \phi (\cdot ,0)\right\vert _{2}^{2}+\left\vert \theta
(\cdot ,0)\right\vert _{2}^{2}\leq C\iint_{Q_{\omega }}\left\vert
\phi \right\vert ^{2}dxdt
\end{equation*}%
for every solution $\left( \phi ,\theta \right) $ of \eqref{adjoint0}.
However, in order to obtain the input space of $L^{\infty }(Q)$, we need to
establish instead an improved observability inequality of the following%
\begin{equation}
\left\vert \phi (\cdot ,0)\right\vert _{2}^{2}+\left\vert \theta (\cdot
,0)\right\vert _{2}^{2}\leq C\iint_{Q_{\omega }}e^{\frac{3}{2}s\alpha
}\left\vert \phi \right\vert ^{2}dxdt,  \label{obs}
\end{equation}%
which can be derived from a global Carleman inequality%
\begin{equation}
\iint_{Q}e^{2s\alpha }\left( \left\vert \phi \right\vert ^{2}+\left\vert
\theta \right\vert ^{2}\right) dxdt\leq C\iint_{Q_{\omega }}e^{\frac{3}{2}%
s\alpha }\left\vert \phi \right\vert ^{2}dxdt \label{cars}
\end{equation}%
for every solution $\left( \phi ,\theta \right) $ of \eqref{adjoint0}. Here,
in \eqref{obs} and \eqref{cars}, $C$ denotes some positive constant
independent of $\phi $ and $\theta $, $\alpha =\alpha (x,t)$ is a weight
function which will be specified precisely in Section 2, and $s$ is a real
number considered as parameter. The basic idea for the inequality %
\eqref{cars} comes originally from \cite{teresa} and \cite{gonzalez}, where
similar inequalities are obtained for some cascaded system and parabolic
system of phase-field.

\bigskip\ Now we state our first result.

\begin{theorem}
\label{th-lin} Let $T>0$. For any $(y_{0},z_{0})\in L^2(\Omega)\times
L^2(\Omega)$, there exists a control $f\in L^{\infty }(Q)$ such that the
solution $(y,z)$ of system \eqref{lineare} corresponding to $f$ satisfies $%
(y,z)\in V^{1}(Q)\times V^1(Q) $ and $y(x,T)=0,\text{\ }z(x,T)=0$ for $x\in
\Omega$ almost everywhere. Moreover, the control $f$ satisfies%
\begin{equation}
\left\Vert f\right\Vert _{\infty }\leq e^{C\kappa }\left( \left\vert
y_{0}\right\vert _2+\left\Vert z_{0}\right\Vert _2\right)
\label{fcontrol}
\end{equation}%
where $C$ is a positive constant depending only on $\Omega $ and $\omega $,
and%
\begin{equation}
\kappa =(1+\Vert a\Vert _{\infty }^{2}+\Vert B\Vert _{\infty }^{2})T+{\frac{1%
}{T}}+1+\Vert a\Vert _{\infty }+\Vert B\Vert _{\infty }.  \label{kconst}
\end{equation}
\end{theorem}

The approach used here to obtain the $L^{\infty }(Q)$ control is originally
from \cite{barbu} (see also \cite{zhang}). We improve this approach to get
the explicit representation of the bound with respect to $T$ by adopting
some techniques from semigroup theory such as $L^{p}$-$L^{q}$ estimate and
maximal $L^{p}$-regularity.

\bigskip The main result of this paper is the following Theorem \ref{th-non}.

\begin{theorem}
\label{th-non} Let $p>N+2$. Let $(\overline{u},\overline{v})$ be a
trajectory of \eqref{ks} corresponding to $(\bar{u}_{0},\bar{v}_{0})$ and
satisfy \eqref{regularity}. Then, there exists a positive constant $c_{1}$
independent of $T$ such that for each $(u_{0},v_{0})$ that satisfies
\begin{equation}
\left\vert u_{0}-\overline{u}_{0}\right\vert _{\infty }+\left\Vert v_{0}-%
\overline{v}_{0}\right\Vert _{W^{2(1-{\frac{1}{p}}),p}(\Omega )}\leq
e^{-c_{1}\left( 1+T+{\frac{1}{T}}\right) },  \label{intialdata}
\end{equation}%
there is a control $f\in L^{\infty }(Q)$ such that system \eqref{e} admits a
solution $(u,v)$ satisfying%
\begin{equation*}
u\in V^{1}(Q)\cap L^{\infty }(Q),v\in V^{2}(Q)\cap L^{\infty }(Q),
\end{equation*}%
and $u(x,T)=\overline{u}(x,T),v(x,T)=\overline{v}(x,T)$ for $x\in \Omega $
almost everywhere.
\end{theorem}

\bigskip

We proceed as follows. In next section, Section 2, we give some preliminary
results. Section 3 is devoted to the proof of the Theorem \ref{th-lin}. The
proof of Theorem \ref{th-non} is presented in section 4.

It is pointed out that throughout the paper, we use $C$ to denote a
positive constant that is independent of time $T$ in most cases but
may be dependent of $\Omega ,\omega $. In the later case we may
write $C(\Omega,\omega)$ instead of a  special specification.

\bigskip

\section{Preliminaries}

In this section, we collect some results that are needed in later sections.
These results are particularly useful in the establishment of the regularity
of linear parabolic system and the $L^\infty$-estimate of controls.

For $p\in (1,\infty ),$ let $A:=A_{p}$ denote the sectorial operator defined
by
\begin{equation}
A_{p}u:=-\Delta u, \; \forall \; u\in D(A_{p}):=\left\{ u\in W^{2,p}(\Omega
);\partial _{\nu}u|_{\partial \Omega }=0\right\} .  \label{laplace}
\end{equation}
Suppose that $\gamma $ is a positive constant.

\begin{itemize}
\item[(i)] Let $\alpha \geq 0$ and $D\left( (A+\gamma )^{\alpha }\right) $
be the function space endowed with the graph norm. Then $D\left( (A+\gamma
)^{\alpha }\right) $ is a Banach space with the following embedding
properties (\cite[p.39]{henry})%
\begin{eqnarray}
D\left( (A+\gamma )^{\alpha }\right) &\hookrightarrow &W^{1,p}(\Omega )\text{
\ if }\alpha >\frac{1}{2},  \label{domain1} \\
\text{and \ }D\left( (A+\gamma )^{\alpha }\right) &\hookrightarrow
&C^{\gamma }(\overline{\Omega })\text{ \ if }0\leq \gamma <2\alpha -\frac{n}{%
p}.  \label{domain2}
\end{eqnarray}

\item[(ii)] Let $\left\{ e^{-tA}\right\} _{t\geq 0}$ \ and $\left\{
e^{-t(A+\gamma )}\right\} _{t\geq 0}$ be the analytic $C_{0}$-semigroups
generated by $-A$ and $-(A+\gamma )$ on $L^{p}(\Omega )(1<p<\infty )$,
respectively. By standard $C_{0}$-semigroup theory, we have (\cite%
{davies,rothe})
\begin{eqnarray}
\left\vert e^{-tA}u\right\vert _{q} &\leq &Cm(t)^{-\frac{N}{2}(\frac{1}{p}-%
\frac{1}{q})}\left\vert u\right\vert _{p},  \label{semigroup1} \\
\text{and \ }\left\vert (A+\gamma )^{\alpha }e^{-t(A+\gamma )}\right\vert
_{q} &\leq &Ct^{-\frac{N}{2}(\frac{1}{p}-\frac{1}{q})-\alpha }\left\vert
u\right\vert _{p}  \label{semigroup2}
\end{eqnarray}%
for all $u\in L^{p}(\Omega ),t>0$ and $1<p\leq q<\infty ,$ where $m(t)=\min
\{1,t\}$.

\item[(iii)] Let $\alpha \geq 0$ and $1<p<\infty $. Then for any $%
\varepsilon >0$ there exists a constant $C_{\varepsilon }$ depending on $%
\Omega $, $\varepsilon $ and $p$ such that (\cite[Lemma 2.1]{horstmann})%
\begin{equation}
\left\vert (A+\gamma )^{\alpha }e^{-tA}\nabla \cdot u\right\vert _{p}\leq
C_{\varepsilon }t^{-\alpha -\frac{1}{2}-\varepsilon }\left\vert u\right\vert
_{p}  \label{semigroup3}
\end{equation}%
for all $u\in L^{p}(\Omega ),t>0$.
\end{itemize}

As a consequence of \eqref{semigroup1} and \eqref{semigroup3}, we
have

\begin{itemize}
\item[(iv)] For any $\varepsilon >0$, there exists a constant $%
C_{\varepsilon }$ depending on $\Omega $, $\varepsilon $ and $p$, such that
\begin{equation}
\left\vert e^{-tA}\nabla \cdot u\right\vert _{q}\leq C_{\varepsilon }m(t)^{-%
\frac{1}{2}-\varepsilon
-\frac{N}{2}(\frac{1}{p}-\frac{1}{q})}\left\vert u\right\vert _{p}
\label{semigroup4}
\end{equation}%
for all $u\in L^{p}(\Omega )$, $t>0$, $1<p\leq q<\infty $.
\end{itemize}

\begin{itemize}
\item[(v)] (\textit{Maximal regularity}) Let $1<p<\infty$. If $F\in L^{p}(Q)$
and $u_{0}\in W^{2(1-{\frac{1}{p}}),p}(\Omega )$ with $\partial _{\nu
}u_{0}=0$ on $\partial \Omega $, then there exists a unique solution of
\begin{equation*}
\frac{du}{dt}=(A+\gamma )u+F\text{ for a.e. }t\in (0,T),\text{ }u(0)=u_{0}
\end{equation*}%
that satisfies
\begin{equation}
\left\Vert \frac{du}{dt}\right\Vert _{p}^{p}+\left\Vert (A+\gamma
)u\right\Vert _{p}^{p}+\left\Vert u\right\Vert _{p}^{p}\leq C\left(
\left\Vert F\right\Vert _{p}^{p}+\left\Vert u_{0}\right\Vert _{W^{2(1-\frac{1%
}{p}),p}(\Omega )}^{p}\right) ,  \label{mr}
\end{equation}%
where $C$ is a positive constant independent of $T$ and $F$.
\end{itemize}

\noindent Inequality (\ref{mr}) was first established as Theorem 9.1 of \cite%
{lady} in Chapter IV, but the independency of $C$ with respect to $T$ is
given later as Theorem 1.1 of \cite{lamberton} (see also Theorem 2.3 of \cite%
{GigaSohr}). \bigskip

Now we consider the well-posedness of the following linear parabolic system
which contains (\ref{lineare}) as its special case.
\begin{equation}
\begin{cases}
\partial _{t}y=\Delta y-\nabla \cdot \left( By\right) -\nabla \cdot \left(
a\nabla z\right) +F & \mathrm{in}\ Q, \\
\partial _{t}z=\Delta z-\gamma z+\delta y & \mathrm{in}\ Q, \\
{\partial }_{\nu }y=0,\partial _{\nu }z=0 & \mathrm{on}\ \Sigma , \\
y(x,0)=y_{0}(x)\ \ z(x,0)=z_{0}(x) & x\in \Omega .%
\end{cases}
\label{linear}
\end{equation}

\begin{proposition}
\label{pro1} Let $a\in L^{\infty }(Q)$ and $B\in L^{\infty }(Q)^{N}$ with $%
B\cdot \nu =0$ on $\Sigma $.

\begin{itemize}
\item[(i)] If $y_{0},z_{0}\in L^{2}(\Omega )$ and $F\in L^{2}(Q)$, then
system \eqref{linear} admits a unique solution $(y,z)\in V^{1}(Q)\times
V^{1}(Q)$ satisfying
\begin{equation}
\Vert y\Vert _{V^{1}(Q)}^{2}+\Vert z\Vert _{V^{1}(Q)}^{2}\leq e^{C\kappa
}\left( \left\vert y_{0}\right\vert _{2}^{2}+\left\vert z_{0}\right\vert
_{2}^{2}+\Vert F\Vert _{2}^{2}\right) ;  \label{prol2}
\end{equation}

\item[(ii)] Let $2\leq p<\infty $. If $F\in L^{p}(Q)$, $y_{0}\in
L^{p}(\Omega )$ and $z_{0}\in W^{2(1-\frac{1}{p}),p}(\Omega )$ with $%
\partial _{\nu }z_{0}=0$ on $\partial \Omega $, then  system %
\eqref{linear} admits a unique solution $(y,z)\in L^{p}(Q)\times
W_{p}^{2,1}(Q)$ satisfying
\begin{equation}
\left\Vert y\right\Vert _{p}^{p}+\left\Vert z\right\Vert
_{W_{p}^{2,1}(Q)}^{p}\leq e^{C\kappa }\left( |y_{0}|_{p}^{p}+\Vert
z_{0}\Vert _{W^{2(1-\frac{1}{p}),p}(\Omega )}^{p}+\left\Vert F\right\Vert
_{p}^{p}\right) ;  \label{prolp}
\end{equation}

\item[(iii)] Let $p>N+2$. If $F\in L^{\infty }(Q)$, $y_{0}\in L^{\infty
}(\Omega )$ and $z_{0}\in W^{1,p}(\Omega )$ with $\partial _{\nu }z_{0}=0$
on $\partial \Omega $, then  system \eqref{linear} admits a solution $%
(y,z)\in L^{\infty }(Q)\times L^{\infty }(Q)$ satisfying
\begin{equation}
\Vert y\Vert _{\infty }+\Vert z\Vert _{\infty }\leq e^{C\kappa }\left(
|y_{0}|_{\infty }+\Vert z_{0}\Vert _{W^{1,p}(\Omega )}+\Vert F\Vert _{\infty
}\right) ,  \label{thlinfty}
\end{equation}
\end{itemize}

\noindent where $\kappa $ is given by \eqref{kconst} and $C=C(\Omega )$.
\end{proposition}

\noindent \textit{Proof.}\quad The existence of solution with respect to $%
y_{0},z_{0}$ and $F$ in different function spaces can be deduced similarly
as in \cite{lady} for which we omit here. We only show the required
estimates with respect to time $T$. Since the proof for \eqref{prol2} is
similar to \eqref{prolp}, we need only to show \eqref{prolp}. Multiply the
first equation of \eqref{linear} by $|y|^{p-2}y$ and integrate over $\Omega $%
, to get%
\begin{equation}
\frac{d}{dt}\left\vert y\right\vert _{p}^{p}+\int_{\Omega }\left\vert \nabla
y\right\vert ^{2}\left\vert y\right\vert ^{p-2}dx\leq C\left( 1+\left\Vert
a\right\Vert _{\infty }^{2}+\left\Vert B\right\Vert _{\infty }^{2}\right)
\left\vert y\right\vert _{p}^{p}+C\left\Vert a\right\Vert _{\infty
}^{2}\left\vert \nabla z\right\vert _{p}^{p}+C\left\vert F\right\vert
_{p}^{p},  \label{th1p0}
\end{equation}%
and in the same way, to get from the second equation of \eqref{linear} that%
\begin{equation}
\frac{d}{dt}\left\vert z\right\vert _{p}^{p}+\int_{\Omega }\left\vert \nabla
z\right\vert ^{2}\left\vert z\right\vert ^{p-2}dx+\left\vert z\right\vert
_{p}^{p}\leq C\left\vert y\right\vert _{p}^{p}.  \label{th1p1}
\end{equation}%
Differentiate $\left\vert \nabla z\right\vert _{p}^{p}$ with respect to $t$
and take the second equation of \eqref{linear} into account again to obtain
\begin{equation}
\frac{d}{dt}\left\vert \nabla z\right\vert _{p}^{p}+\int_{\Omega }\left\vert
\nabla z\right\vert ^{p-2}\left\vert \Delta z\right\vert ^{2}dx\leq
C\left\vert \nabla z\right\vert _{p}^{p}+C\left( \left\vert y\right\vert
_{p}^{p}+\left\vert z\right\vert _{p}^{p}\right) .  \label{thlp2}
\end{equation}%
The inequalities \eqref{th1p0}-\eqref{thlp2} together with Gronwall's
inequality lead to
\begin{equation}
\left\vert y(\cdot ,t)\right\vert _{p}^{p}+\left\vert z(\cdot ,t)\right\vert
_{p}^{p}+\left\vert \nabla z(\cdot ,t)\right\vert _{p}^{p}\leq e^{C\kappa
}\left( \left\vert y_{0}\right\vert _{p}^{p}+\left\Vert z_{0}\right\Vert
_{W^{1,p}(\Omega )}^{p}\right)  \label{th1pp}
\end{equation}%
for all $t\in \lbrack 0,T]$. On the other hand, by the maximal regularity %
\eqref{mr} for the second equation of \eqref{linear}, it follows that
\begin{equation*}
\left\Vert \partial _{t}z\right\Vert _{p}^{p}+\left\Vert \Delta z\right\Vert
_{p}^{p}+\left\Vert z\right\Vert _{p}^{p}\leq C\left( \left\Vert
z_{0}\right\Vert _{W^{2(1-\frac{1}{p}),p}(\Omega )}^{p}+\left\Vert
y\right\Vert _{p}^{p}+\left\Vert z\right\Vert _{p}^{p}\right) ,
\end{equation*}%
which together with \eqref{th1pp} yields \eqref{prolp}.

Now we turn to the L$^{\infty }$-estimate \eqref{thlinfty}. We first assume
that $y_{0}\in C\left( \overline{\Omega }\right) $ and $F\in C\left(
\overline{Q}\right) $. Let $A$ be defined by \eqref{laplace}, and let $%
\left\{ e^{-tA}\right\} _{t\geq 0}$ and $\left\{ e^{-t(A+\gamma )}\right\}
_{t\geq 0}$ be the analytic $C_{0}$-semigroups generated by $-A$ and $%
-(A+\gamma )$ in $L^{p}(\Omega ),1<p<\infty $, respectively. Then the
solution $(y,z)$ of system \eqref{linear} can be represented as follows%
\begin{eqnarray}
y(\cdot ,t) &=&e^{-tA}y_{0}+\int_{0}^{t}e^{-(t-s)A}\left[ -\nabla \cdot
(By)-\nabla \cdot \left( a\nabla z\right) +F\right] (\cdot ,s)ds,
\label{ysol0} \\
z(\cdot ,t) &=&e^{-t(A+\gamma )}z_{0}+\delta \int_{0}^{t}e^{-(t-s)(A+\gamma
)}y(\cdot ,s)ds.  \label{zsol0}
\end{eqnarray}%
Take the norm of $C\left( \overline{\Omega }\right) $ on both sides of %
\eqref{ysol0} to get
\begin{eqnarray}
\left\Vert y(\cdot ,t)\right\Vert _{C\left( \overline{\Omega }\right) }
&\leq &\left\Vert e^{-tA}y_{0}\right\Vert _{C\left( \overline{\Omega }%
\right) }+\int_{0}^{t}\left\Vert e^{-(t-s)A}\nabla \cdot \left( By+a\nabla
z\right) (\cdot ,s)\right\Vert _{C\left( \overline{\Omega }\right) }ds
\notag \\
&&+\int_{0}^{t}\left\Vert e^{-(t-s)A}F(\cdot ,s)\right\Vert _{C\left(
\overline{\Omega }\right) }ds.  \label{th1p3}
\end{eqnarray}%
To estimate \eqref{th1p3}, we first observe that the operator $-A$ generates
a bounded analytic semigroup on $C\left( \overline{\Omega }\right) $(\cite%
{arendt}). It follows from the maximum principle that
\begin{equation}
\left\Vert e^{-tA}y_{0}\right\Vert _{C\left( \overline{\Omega }\right) }\leq
\left\Vert y_{0}\right\Vert _{C\left( \overline{\Omega }\right) },
\label{th1p4}
\end{equation}%
and%
\begin{equation}
\left\Vert e^{-(t-s)A}F(\cdot ,s)\right\Vert _{C\left( \overline{\Omega }%
\right) }\leq \left\Vert F(\cdot ,s)\right\Vert _{C\left( \overline{\Omega }%
\right) }  \label{thlp5}
\end{equation}%
for any $0\leq s\leq t.$ Since $p>N+2$, we can take $\varepsilon $ and $%
\alpha $ such that
\begin{equation*}
0<\varepsilon <\frac{p-N-2}{2p}\text{ \ and \ }\frac{N}{2p}<\alpha <\frac{1}{%
2}-\frac{1}{p}-\varepsilon .
\end{equation*}%
Then, with the help of \eqref{domain2}, \eqref{semigroup3} and the
H\"{o}lder inequality, we have, for any $t\in \lbrack 0,T],$ that
\begin{eqnarray*}
&&\int_{0}^{t}\left\Vert e^{-(t-s)A}\nabla \cdot \left( By+a\nabla z)\right)
(\cdot ,s)\right\Vert _{C\left( \overline{\Omega }\right) }ds \\
&\leq &\int_{0}^{t}\left\vert (A+\gamma )^{\alpha }e^{-(t-s)A}\left(
By+a\nabla z\right) (\cdot ,s)\right\vert _{p}ds \\
&\leq &C\int_{0}^{t}(t-s)^{-\alpha -\frac{1}{2}-\varepsilon }\left\vert
\left( By+a\nabla z\right) (\cdot ,s)\right\vert _{p}ds \\
&\leq &C\left( 1+\left\Vert a\right\Vert _{\infty }+\left\Vert B\right\Vert
_{\infty }\right) (\left\Vert y\right\Vert _{p}+\left\Vert \nabla
z\right\Vert _{p})T^{\frac{1}{2}-\alpha -\varepsilon -\frac{1}{p}}.
\end{eqnarray*}%
This together with \eqref{th1pp} gives
\begin{equation}
\int_{0}^{t}\left\Vert e^{-(t-s)A}\nabla \cdot \left( By+a\nabla z)\right)
(\cdot ,s)\right\Vert _{C\left( \overline{\Omega }\right) }ds\leq e^{C\kappa
}\left( \left\vert y_{0}\right\vert _{p}+\left\Vert z_{0}\right\Vert
_{W^{1,p}(\Omega )}+\left\Vert F\right\Vert _{p}\right) .  \label{th1p6}
\end{equation}%
By \eqref{th1p3}-\eqref{th1p6}, we obtain
\begin{equation}
\left\Vert y\right\Vert _{\infty }\leq e^{C\kappa }\left( \left\Vert
y_{0}\right\Vert _{\infty }+\left\Vert z_{0}\right\Vert _{W^{1,p}(\Omega
)}+\left\Vert F\right\Vert _{\infty }\right) .  \label{th1p7}
\end{equation}%
Next, take the norm of $W^{1,p}(\Omega )$ on both sides of \eqref{zsol0} to
get
\begin{eqnarray}
\left\Vert z(\cdot ,t)\right\Vert _{W^{1,p}(\Omega )} &\leq &\left\Vert
e^{-t(A+\gamma )}z_{0}\right\Vert _{W^{1,p}(\Omega )}  \notag \\
&&+\delta \int_{0}^{t}\left\Vert e^{-(t-s)(A+\gamma )}y(\cdot ,s)\right\Vert
_{W^{1,p}(\Omega )}ds,  \label{thlp7b}
\end{eqnarray}%
for any $t\in \lbrack 0,T].$ To estimate \eqref{thlp7b}, we first notice
that
\begin{equation}
\left\Vert e^{-t(A+\gamma )}z_{0}\right\Vert _{W^{1,p}(\Omega )}\leq
e^{CT}\left\Vert z_{0}\right\Vert _{W^{1,p}(\Omega )}  \label{thlp8}
\end{equation}%
which can be obtained by the same energy method used in proving \eqref{th1pp}%
. Let $\frac{1}{2}<\alpha <1-\frac{1}{p}.$ By \eqref{domain1}, %
\eqref{semigroup2} and the H\"{o}lder inequality, we have that for any $t\in
\lbrack 0,T],$%
\begin{eqnarray*}
\int_{0}^{t}\left\Vert e^{-(t-s)(A+\gamma )}y(\cdot ,s)\right\Vert
_{W^{1,p}(\Omega )}ds &\leq &C\int_{0}^{t}\left\vert (A+\gamma )^{\alpha
}e^{-(t-s)(A+\gamma )}y(\cdot ,s)\right\vert _{p}ds \\
&\leq &C\int_{0}^{t}(t-s)^{-\alpha }\left\vert y(\cdot ,s)\right\vert _{p}ds
\\
&\leq &C\left\Vert y\right\Vert _{p}T^{-\alpha +1-\frac{1}{p}}.
\end{eqnarray*}
This together with \eqref{th1pp} gives
\begin{equation}
\int_{0}^{t}\left\Vert e^{-(t-s)(A+\gamma )}y(\cdot ,s)\right\Vert
_{W^{1,p}(\Omega )}ds\leq e^{C\kappa }\left( \left\vert y_{0}\right\vert
_{p}+\left\Vert z_{0}\right\Vert _{W^{1,p}(\Omega )}+\left\Vert F\right\Vert
_{p}\right) .  \label{thlp9}
\end{equation}%
Finally, by \eqref{thlp7b}-\eqref{thlp9} and the Sobolev embedding $%
W^{1,p}(\Omega )\hookrightarrow C(\overline{\Omega })$ for $p>N$, we get%
\begin{equation}
\left\Vert z\right\Vert _{\infty }\leq e^{C\kappa }\left( \left\Vert
y_{0}\right\Vert _{\infty }+\left\Vert z_{0}\right\Vert _{W^{1,p}(\Omega
)}+\left\Vert F\right\Vert _{\infty }\right) .  \label{thlp10}
\end{equation}

To complete the proof, let us consider the general case that
$y_{0}\in L^{\infty }(\Omega )$ and $F\in L^{\infty }(Q)$. This can
be done by smoothing the data and density  argument.
Precisely, let $\left\{ y_{0n}\right\} _{n=1}^{\infty }\subset C(\bar{\Omega}%
)$ and $\left\{ F_{n}\right\} _{n=1}^{\infty }\subset C(\bar{Q})$ be such
that $y_{0n}\rightarrow y_{0}$ in $L^{2}(\Omega )$, $F_{n}\rightarrow F$ in $%
L^{2}(Q)$ and $\left\vert y_{0n}\right\vert _{\infty }\leq \left\vert
y_{0}\right\vert _{\infty },\left\Vert F_{n}\right\Vert _{\infty }\leq
\left\Vert F\right\Vert _{\infty }.$ For each $n$, let $(y_{n},z_{n})$ be a
solution of \eqref{linear} corresponding to $y_{0n},z_{0},F_{n}$, which
satisfies the inequalities \eqref{prol2} and \eqref{thlinfty} with $(y,z)$
replaced by $(y_{n},z_{n})$. Thus, by the uniformly boundedness, we can
extract subsequences of $(y_{n},z_{n})$ such that it converges to $(y,z),$
which is a weak solution of \eqref{linear} corresponding to $y_{0},z_{0}$
and $F$. Moreover, $y,z\ $satisfy the inequality \eqref{thlinfty}. \hfill $%
\Box $

\section{ Proof of Theorem \protect\ref{th-lin}}

To prove Theorem \ref{th-lin}, we first establish a global Carleman
inequality for the adjoint system \eqref{adjoint0}.

Let $\omega ^{\prime }\subset \subset \omega $, that is, $\overline{\omega
^{\prime }}\subset \omega $, be a nonempty open subset. Then, there is a
function $\beta \in C^{2}(\overline{\Omega })$ such that $\beta (x)>0$ for
all $x\in \Omega $, and $\beta |_{\partial \Omega }=0,\left\vert \nabla
\beta (x)\right\vert >0$ for all $x\in \overline{\Omega \setminus \omega
^{\prime }}$ (see \cite[Lemma 1.1]{fursikov}). For $\lambda >0,$  set%
\begin{equation}
\varphi =\frac{e^{\lambda \beta }}{t(T-t)},\text{ \ }\alpha =\frac{%
e^{\lambda \beta }-e^{2\lambda \left\Vert \beta \right\Vert _{C(\overline{%
\Omega })}}}{t(T-t)},  \label{alpha}
\end{equation}%
and
\begin{equation}
\gamma (\lambda )=e^{2\lambda \left\Vert \beta \right\Vert _{C(\overline{%
\Omega })}}.  \label{gamma}
\end{equation}

\begin{lemma}
\label{carleman} Let $f_{i}\in L^{2}(Q)$, $i=0,1,\ldots ,N$. Then there
exists a constant $\lambda _{0}=\lambda _{0}(\Omega ,\omega ^{\prime })>1$,
such that for all $\lambda \geq \lambda _{0}$ and $s\geq \gamma (\lambda
)(T+T^{2})$,
\begin{eqnarray}
&&\iint_{Q}\left[ (s\varphi )^{1+d}|\nabla z|^{2}+(s\varphi )^{3+d}|z|^{2}%
\right] e^{2s\alpha }\ dxdt  \notag \\
&\leq &C\left(\iint_{Q}(s\varphi )^{d}e^{2s\alpha }|f_{0}|^{2}\
dxdt+\sum_{i=1}^{N}\iint_{Q}(s\varphi )^{2+d}e^{2s\alpha
}|f_{i}|^{2}\ dxdt \notag \right.\\
&&+\left.\iint_{Q_{\omega ^{\prime }}}(s\varphi )^{3+d}e^{2s\alpha }|z|^{2}\ dxdt%
\right) \label{carle0}
\end{eqnarray}%
for all solutions $z$ to the equation
\begin{equation*}
\begin{cases}
\partial _{t}z-\Delta z=f_{0}+\sum_{i=1}^{N}{\frac{\partial f_{i}}{\partial
x_{i}}} & \mathrm{in}\ Q, \\
\partial _{\nu }z=0\ \  & \mathrm{on}\ \Sigma , \\
z(x,0)=z_{0}(x) & x\in \Omega ,%
\end{cases}%
\end{equation*}%
with $z_{0}\in L^{2}(\Omega )$, where $C=C(\Omega ,\omega ^{\prime })$, and $%
\gamma (\lambda )$ given by \eqref{gamma}.
\end{lemma}

Essentially speaking, Lemma \ref{carleman} has been proven in \cite%
{imanuvilov} (see also \cite{gonzalez}) but the explicit independency of the
constant $C$ with respect to $T$ is shown in a similar way as in \cite%
{fernandez1} and \cite{fernandez4}. For notational simplicity in the sequel,
we introduce
\begin{equation}
I_{1}(s,\lambda ;\phi )=\iint_{Q}\left[ (s\varphi )^{3}\left\vert \nabla
\phi \right\vert ^{2}+(s\varphi )^{5}\left\vert \phi \right\vert ^{2}\right]
e^{2s\alpha }dxdt,  \label{i1}
\end{equation}%
and%
\begin{equation}
I_{2}(s,\lambda ;\theta )=\iint_{Q}\left[ s\varphi \left\vert \nabla \theta
\right\vert ^{2}+(s\varphi )^{3}\left\vert \theta \right\vert ^{2}\right]
e^{2s\alpha }dxdt.  \label{i2}
\end{equation}

\begin{lemma}
\label{carleman1} There exists a positive constant $\lambda _{1}=C(\Omega
,\omega ,\omega ^{\prime })(1+\left\Vert a\right\Vert _{\infty
}^{2}+\left\Vert B\right\Vert _{\infty }^{2})$ satisfying $\gamma (\lambda
_{1})\geq \lambda _{1}>1$ such that for any $\lambda \geq \lambda _{1},s\geq
\gamma (\lambda )(T+T^{2})$ and $\phi ^{T},\theta ^{T}\in L^{2}(\Omega ),$
the associated solution $(\phi ,\theta )$ to \eqref{adjoint0} satisfies%
\begin{equation}
I_{1}(s,\lambda ;\phi )+I_{2}(s,\lambda ;\theta )\leq C_{1}\iint_{Q_{\omega
}}\lambda ^{8}(s\varphi )^{9}e^{2s\alpha }\left\vert \phi \right\vert
^{2}dxdt,  \label{carle1}
\end{equation}%
where $C_{1}=C_{1}(\Omega ,\omega ^{\prime },\omega ).$
\end{lemma}

\noindent\textit{Proof.}\quad Applying Lemma \ref{carleman} to the first
equation of \eqref{adjoint0} with $d=2$ and the second one with $d=0$,
respectively, we obtain that there exist positive constants $%
c_{0}(\Omega,\omega ^{\prime })$ and $\lambda _{1}^{0}$ satisfying%
\begin{equation}
\gamma (\lambda _{1}^{0})\geq \lambda _{1}^{0}=c_{0}(\Omega ,\omega ^{\prime
})\left( 1+\left\Vert a\right\Vert _{\infty }^{2}+\left\Vert B\right\Vert
_{\infty }^{2}\right) >1  \label{lambda1}
\end{equation}%
such that for all $\lambda \geq \lambda _{1}^{0}$ and $s\geq \gamma \left(
\lambda \right) (T+T^{2}),$%
\begin{equation}
I_{1}(s,\lambda ;\phi )+I_{2}(s,\lambda ;\theta )\leq c_{1}\iint_{Q_{\omega
^{\prime }}}\left[ (s\varphi )^{5}\left\vert \phi \right\vert ^{2}+(s\varphi
)^{3}\left\vert \theta \right\vert ^{2}\right] e^{2s\alpha }dxdt
\label{i1i2a}
\end{equation}%
for all solutions $(\phi ,\theta )$ to \eqref{adjoint0} with $\phi
^{T},\theta ^{T}\in L^{2}(\Omega ),$ where and in what follows, the symbol $%
c_{i},i=1,2,\ldots,$ stand for some positive constants depending on $%
\Omega,\omega ^{\prime }$ and $\omega $.

Next, let $\xi \in C_{0}^{\infty }(\Omega )$ be such that $\xi =1$ in $%
\omega ^{\prime }$, $\xi =0$ in $\Omega \setminus \overline{\omega }$, $%
0\leq \xi \leq 1$ in $\omega $, and
\begin{equation}
\Delta \xi \cdot \xi ^{-1/2}\in L^{\infty }(\Omega ),\nabla \xi \cdot \xi
^{-1/2}\in L^{\infty }(\Omega )^{N}.  \label{xi}
\end{equation}%
The existence of such a function $\xi $ is easy to obtain (see, for instance
\cite{teresa}). Set%
\begin{equation*}
\eta =(s\varphi )^{3}e^{2s\alpha }.
\end{equation*}%
Multiply the first equation of \eqref{adjoint0} by $\theta \eta \xi $ to get
\begin{eqnarray*}
&&\delta \iint_{Q}(s\varphi )^{3}e^{2s\alpha }\left\vert \theta \right\vert
^{2}\xi dxdt=\iint_{Q}\eta \xi \theta \left[ -\partial _{t}\phi -\Delta \phi
-B\nabla \phi \right] dxdt \\
&=&\iint_{Q}\left\{ \eta \xi \phi \left[ -\Delta \theta +\gamma \theta
+\nabla \cdot \left( a\nabla \phi \right) \right] +\phi \theta \xi \left(
\partial _{t}\eta \right) +\eta \xi \theta \left( -\Delta \phi -B\nabla \phi
\right) \right\} dxdt.
\end{eqnarray*}%
Integration by parts gives%
\begin{equation}
\delta \iint_{Q}(s\varphi )^{3}e^{2s\alpha }\left\vert \theta \right\vert
^{2}\xi dxdt=\sum_{i=1}^{7}J_{i},  \label{jsum}
\end{equation}%
where
\begin{equation*}
\begin{array}{ll}
\displaystyle J_{1}=\iint_{Q}\phi \theta (\partial _{t}\eta +\gamma \eta
)\xi dxdt, & \displaystyle J_{2}=\iint_{Q}\phi \nabla (\eta \xi )\cdot
\nabla \theta dxdt, \\
\displaystyle J_{3}=-\iint_{Q}a\phi \nabla (\eta \xi )\cdot \nabla \phi dxdt,
& \displaystyle J_{4}=\iint_{Q}\theta \nabla (\eta \xi )\cdot \nabla \phi
dxdt, \\
\displaystyle J_{5}=-\iint_{Q}\theta \eta \xi B\nabla \phi dxdt, & %
\displaystyle J_{6}=\iint_{Q}\eta \xi \nabla \theta \cdot \nabla \phi dxdt,
\\
\displaystyle J_{7}=-\iint_{Q}a\left\vert \nabla \phi \right\vert ^{2}\eta
\xi dxdt. &
\end{array}%
\end{equation*}%
To estimate these integrals, we first observe by \eqref{alpha} and \eqref{xi}
that%
\begin{equation*}
\left\vert \partial _{t}\eta \right\vert \leq (s\varphi )^{5}e^{2s\alpha };%
\text{ \ }\left\vert \nabla (\eta \xi )\right\vert \leq C(\xi
^{1/2}s^{3}\varphi ^{3}+\xi \lambda s^{4}\varphi ^{4})e^{2s\alpha }.
\end{equation*}%
This together with Cauchy's inequality gives the estimation of $%
J_{i},i=1,\ldots ,6$ as follows:
\begin{equation}
J_{1}\leq \varepsilon _{1}I_{2}(s,\lambda ;\theta )+\frac{C}{4\varepsilon
_{1}}\iint_{Q}\left[ (s\varphi )^{3}+(s\varphi )^{7}\right] e^{2s\alpha
}\left\vert \phi \right\vert ^{2}\xi dxdt;  \label{j1}
\end{equation}%
\begin{equation}
J_{2}\leq \varepsilon _{1}I_{2}(s,\lambda ;\theta )+\frac{C}{4\varepsilon
_{1}}\iint_{Q}\left[ (s\varphi )^{5}+\lambda ^{2}(s\varphi )^{7}\right]
e^{2s\alpha }\left\vert \phi \right\vert ^{2}\xi dxdt;  \label{j2}
\end{equation}%
\begin{equation}
J_{3}\leq \varepsilon _{1}I_{1}(s,\lambda ;\phi )+\frac{C\left\Vert
a\right\Vert _{\infty }^{2}}{4\varepsilon _{1}}\iint_{Q}\left[ (s\varphi
)^{3}+\lambda ^{2}(s\varphi )^{5}\right] e^{2s\alpha }\left\vert \phi
\right\vert ^{2}\xi dxdt;  \label{j3}
\end{equation}%
\begin{equation}
J_{4}\leq \varepsilon _{1}I_{2}(s,\lambda ;\theta )+\frac{c_{2}}{%
2\varepsilon _{1}}\iint_{Q}\left[ (s\varphi )^{3}+\lambda (s\varphi )^{5}%
\right] e^{2s\alpha }\left\vert \nabla \phi \right\vert ^{2}\xi dxdt;
\label{j4}
\end{equation}%
\begin{equation}
J_{5}\leq \varepsilon _{1}I_{2}(s,\lambda ;\theta )+\frac{c_{3}\left\Vert
B\right\Vert _{\infty }^{2}}{4\varepsilon _{1}}\iint_{Q}(s\varphi
)^{3}e^{2s\alpha }\left\vert \nabla \phi \right\vert ^{2}\xi dxdt;
\label{j5}
\end{equation}%
\begin{equation}
J_{6}\leq \varepsilon _{1}I_{2}(s,\lambda ;\theta )+\frac{c_{4}}{%
4\varepsilon _{1}}\iint_{Q}(s\varphi )^{5}e^{2s\alpha }\left\vert
\nabla \phi \right\vert ^{2}\xi dxdt;  \label{j6}
\end{equation}%
\begin{equation}
J_{7}\leq \left\Vert a\right\Vert _{\infty }\iint_{Q}(s\varphi
)^{3}e^{2s\alpha }\left\vert \nabla \phi \right\vert ^{2}\xi dxdt,
\label{j71}
\end{equation}%
where $\varepsilon _{1}$ is an arbitrary positive constant which will be
determined later.

The inequalities \eqref{j4}-\eqref{j71} lead to%
\begin{eqnarray}
J_{4}+J_{5}+J_{6}+J_{7} &\leq &3\varepsilon _{1}I_{2}(s,\lambda ;\theta )+%
\frac{c_{5}}{\varepsilon _{1}}(1+\left\Vert a\right\Vert _{\infty
}^{2}+\left\Vert B\right\Vert _{\infty }^{2})  \notag \\
&&\times \iint_{Q}\lambda ^{2}(s\varphi )^{5}e^{2s\alpha }\left\vert \nabla
\phi \right\vert ^{2}\xi dxdt.  \label{j456}
\end{eqnarray}%
Next, we estimate the integral on the right hand side of the inequality %
\eqref{j456}. Let
\begin{equation*}
\tilde{\eta}=\lambda ^{2}(s\varphi )^{5}e^{2s\alpha }.
\end{equation*}%
Multiply the first equation of \eqref{adjoint0} by $\tilde{\eta}\xi \phi $
and integrate over $Q$ to obtain, by the integration by parts, that
\begin{equation*}
\iint_{Q}\lambda ^{2}(s\varphi )^{5}e^{2s\alpha }\left\vert \nabla \phi
\right\vert ^{2}\xi dxdt=\sum_{i=8}^{11}J_{i},
\end{equation*}%
where
\begin{equation*}
\begin{array}{ll}
\displaystyle J_{8}=-\frac{1}{2}\iint_{Q}\phi ^{2}\xi \partial _{t}\tilde{%
\eta}dxdt, & \displaystyle J_{9}=-\iint_{Q}\phi \nabla (\tilde{\eta}\xi
)\cdot \nabla \phi dxdt, \\
\displaystyle J_{10}=\iint_{Q}\phi \tilde{\eta}\xi B\nabla \phi dxdt, & %
\displaystyle J_{11}=\delta \iint_{Q}\theta \phi \tilde{\eta}\xi dxdt.%
\end{array}%
\end{equation*}%
Since%
\begin{equation*}
\left\vert \partial _{t}\tilde{\eta}\right\vert \leq C\lambda ^{2}(s\varphi
)^{7}e^{2s\alpha },\;\left\vert \nabla (\tilde{\eta}\xi )\right\vert \leq
C(\xi ^{1/2}\lambda ^{2}s^{5}\varphi ^{5}+\lambda ^{3}s^{6}\varphi ^{6}\xi
)e^{2s\alpha },
\end{equation*}%
in the same way of estimating $J_{1}$-$J_{7}$, we can get for any $%
\varepsilon _{2}>0$ that
\begin{equation}
J_{8}\leq C\iint_{Q}\lambda ^{2}(s\varphi )^{7}e^{2s\alpha }\left\vert \phi
\right\vert ^{2}\xi dxdt;  \label{j7}
\end{equation}%
\begin{equation}
J_{9}\leq \varepsilon _{2}I_{1}(s,\lambda ;\phi )+\frac{C}{2\varepsilon _{2}}%
\iint_{Q}\left[ \lambda ^{4}(s\varphi )^{7}+\lambda ^{6}(s\varphi )^{9}%
\right] e^{2s\alpha }\left\vert \phi \right\vert ^{2}\xi dxdt;  \label{j8}
\end{equation}%
\begin{equation}
J_{10}\leq \varepsilon _{2}I_{2}(s,\lambda ;\theta )+\frac{C\left\Vert
B\right\Vert _{\infty }^{2}}{4\varepsilon _{2}}\iint_{Q}\lambda
^{4}(s\varphi )^{7}e^{2s\alpha }\left\vert \phi \right\vert ^{2}\xi dxdt;
\label{j9}
\end{equation}%
\begin{equation}
J_{11}\leq \varepsilon _{2}I_{1}(s,\lambda ;\phi )+\frac{C}{4\varepsilon _{2}%
}\iint_{Q}\lambda ^{4}(s\varphi )^{7}e^{2s\alpha }\left\vert \phi
\right\vert ^{2}\xi dxdt.  \label{j10}
\end{equation}%
Finally, we take
\begin{equation*}
\varepsilon _{1}=\frac{\delta }{10c_{1}},\text{ and }\varepsilon _{2}=\frac{%
\delta }{10c_{1}c_{5}(1+\left\Vert a\right\Vert _{\infty }^{2}+\left\Vert
B\right\Vert _{\infty }^{2})}\times \frac{\delta }{20c_{1}}
\end{equation*}%
to get,  from \eqref{jsum}-\eqref{j10},  that%
\begin{equation*}
I_{1}(s,\lambda ;\phi )+I_{2}(s,\lambda ;\theta )\leq c_{6}(1+\left\Vert
a\right\Vert _{\infty }^{2}+\left\Vert B\right\Vert _{\infty
}^{2})^{2}\iint_{Q_{\omega }}\lambda ^{6}(s\varphi )^{9}e^{2s\alpha
}\left\vert \phi \right\vert ^{2}dxdt.
\end{equation*}%
Thus there is a positive constant%
\begin{equation*}
\gamma (\lambda _{1})\geq \lambda _{1}=c_{6}\left( 1+\left\Vert
a\right\Vert _{\infty }^{2}+\left\Vert B\right\Vert _{\infty
}^{2}\right) \geq \lambda _{1}^{0}>1
\end{equation*}%
such that for any $\lambda \geq \lambda _{1}$ and $s\geq \gamma (\lambda
)(T+T^{2})$, the inequality \eqref{carle1} holds, where $\lambda _{1}^{0}$
is given by \eqref{lambda1}. \hfill $\Box $

\begin{proposition}
\label{pro-obser} There exist positive constants $\lambda $ and $s$ such
that, for all $T>0,\phi ^{T},\theta ^{T}\in L^{2}(\Omega ),$ the solution ($%
\phi ,\theta $) of the system \eqref{adjoint0} satisfies%
\begin{equation}
\left| \phi(\cdot ,0)\right| _{2}^{2}+ \left| \theta (\cdot ,0)\right|
_{2}^{2}\leq e^{C\kappa }\iint_{Q_{\omega }}e^{\frac{3}{2}s\alpha
}\left\vert \phi \right\vert ^{2}dxdt,  \label{observability}
\end{equation}%
where $\kappa $ is given by \eqref{kconst}.
\end{proposition}

\noindent\textit{Proof.}\quad By integration by parts, we observe that%
\begin{equation}
-\frac{d}{dt}\left\vert \phi \right\vert _{2}^{2}+\left\vert \nabla \phi
\right\vert _{2}^{2}\leq (1+\left\Vert B\right\Vert _{\infty
}^{2})\left\vert \phi \right\vert _{2}^{2}+\delta ^{2}\left\vert \theta
\right\vert _{2}^{2},  \label{aj1}
\end{equation}%
and
\begin{equation}
-\frac{d}{dt}\left\vert \theta \right\vert _{2}^{2}+\left\vert \nabla \theta
\right\vert _{2}^{2}+2\gamma \left\vert \theta \right\vert _{2}^{2}\leq
\left\Vert a\right\Vert _{\infty }^{2}\left\vert \nabla \phi \right\vert
_{2}^{2}.  \label{aj2}
\end{equation}%
Suppose first that $\left\Vert a\right\Vert _{\infty }\geq 1$. Multiply %
\eqref{aj1} by $\left\Vert a\right\Vert _{\infty }^{2}$ to get by \eqref{aj2}
that
\begin{equation*}
\frac{d}{dt}\left[ e^{C(1+\left\Vert a\right\Vert _{\infty }^{2}+\left\Vert
B\right\Vert _{\infty }^{2})t}\left( \left\Vert a\right\Vert _{\infty
}^{2}\left\vert \phi \right\vert _{2}^{2}+\left\vert \theta \right\vert
_{2}^{2}\right) \right] \geq 0.
\end{equation*}%
Integrating above inequality over $[0,t]$ for any $t\in (0,T]$ gives
\begin{equation*}
\left\Vert a\right\Vert _{\infty }^{2}\left\vert \phi (\cdot ,0)\right\vert
_{2}^{2}+\left\vert \theta (\cdot ,0)\right\vert _{2}^{2}\leq
e^{C(1+\left\Vert a\right\Vert _{\infty }^{2}+\left\Vert B\right\Vert
_{\infty }^{2})T}\left( \left\Vert a\right\Vert _{\infty }^{2}\left\vert
\phi (\cdot ,t)\right\vert _{2}^{2}+\left\vert \theta (\cdot ,t)\right\vert
_{2}^{2}\right) ,
\end{equation*}%
which implies that
\begin{equation}
\left\vert \phi (\cdot ,0)\right\vert _{2}^{2}+\left\vert \theta (\cdot
,0)\right\vert _{2}^{2}\leq e^{C\left[ (1+\left\Vert a\right\Vert _{\infty
}^{2}+\left\Vert B\right\Vert _{\infty }^{2})T+\left\Vert a\right\Vert
_{\infty }\right] }\left( \left\vert \phi (\cdot ,t)\right\vert
_{2}^{2}+\left\vert \theta (\cdot ,t)\right\vert _{2}^{2}\right)  \label{aj3}
\end{equation}%
for any $t\in (0,T]$. The integration of \eqref{aj3} on both sides over $%
\left[ T/4,3T/4\right] $ leads to
\begin{equation*}
\left\vert \phi (\cdot ,0)\right\vert _{2}^{2}+\left\vert \theta (\cdot
,0)\right\vert _{2}^{2}\leq \frac{2}{T}e^{C\left[ (1+\left\Vert a\right\Vert
_{\infty }^{2}+\left\Vert B\right\Vert _{\infty }^{2})T+\left\Vert
a\right\Vert _{\infty }\right] }\int_{\frac{T}{4}}^{\frac{3T}{4}%
}\int_{\Omega }\left\vert \phi \right\vert ^{2}+\left\vert \theta
\right\vert ^{2}dxdt.
\end{equation*}%
Since
\begin{equation*}
(s\varphi )^{-5}e^{-2s\alpha },(s\varphi )^{-3}e^{-2s\alpha }\leq e^{\frac{Cs%
}{T^{2}}}\text{ in }\Omega \times \left[ \frac{T}{4},\frac{3T}{4}\right] ,
\end{equation*}%
it follows by \eqref{carle1} that
\begin{equation*}
\left\vert \phi (\cdot ,0)\right\vert _{2}^{2}+\left\vert \theta (\cdot
,0)\right\vert _{2}^{2}\leq \frac{2C_{1}}{T}e^{C\left[ (1+\left\Vert
a\right\Vert _{\infty }^{2}+\left\Vert B\right\Vert _{\infty
}^{2})T+\left\Vert a\right\Vert _{\infty }\right] +\frac{Cs}{T^{2}}%
}\iint_{Q_{\omega }}\lambda ^{8}(s\varphi )^{9}e^{2s\alpha }\left\vert \phi
\right\vert ^{2}dxdt,
\end{equation*}%
where by taking $\lambda $ and $s$ as
\begin{equation*}
\lambda =C\left( 1+\left\Vert a\right\Vert _{\infty }^{2}+\left\Vert
B\right\Vert _{\infty }^{2}\right) ,s=C\left( 1+\left\Vert a\right\Vert
_{\infty }^{2}+\left\Vert B\right\Vert _{\infty }^{2}\right) (T+T^{2}),
\end{equation*}%
we get \eqref{observability}.

Finally, if $\Vert a\Vert _{\infty }<1$, then
\begin{equation*}
\frac{d}{dt}\left[ e^{C(1+\left\Vert B\right\Vert _{\infty }^{2})t}\left(
\left\vert \phi \right\vert _{2}^{2}+\left\vert \theta \right\vert
_{2}^{2}\right) \right] \geq 0.
\end{equation*}%
is a direct consequence of \eqref{aj1} and \eqref{aj2}. Thus, \eqref{aj3} verifies. In a similar
argument as in the proof of $\Vert a\Vert _{\infty }\ge 1$, one can
easily get \eqref{observability}. This completes the proof.\hfill $%
\Box $

\vspace{0.3cm}

\noindent\textbf{Proof of Theorem \ref{th-lin}.} Let $s$ and $\lambda $ be
such that the observability estimate \eqref{observability} and
\begin{equation}
\eta (\lambda )=e^{-\lambda \left\Vert \beta \right\Vert _{C(\overline{\Omega%
})}}<\frac{1}{2}  \label{eta}
\end{equation}%
hold. Let $\varepsilon >0$ and consider the following optimal control problem%
\begin{equation*}
\text{Minimize}\left\{ \iint_{Q_{\omega }}\left\vert f\right\vert ^{2}e^{-%
\frac{3}{2}s\alpha }dxdt+\frac{1}{\varepsilon }\left( \left\vert y(\cdot
,T)\right\vert _{2}^{2}+\left\vert z(\cdot ,T)\right\vert _{2}^{2}\right)
\right\}
\end{equation*}%
subject to all $f\in L^{2}(Q)$, where $(y,z)$ is the solution of %
\eqref{lineare} associated to $f$. The existence of an optimal pair $%
(f_{\varepsilon },y_{\varepsilon },z_{\varepsilon })$ to the above optimal
control problem follows from the standard argument. By the Pontryagin
maximum principle (\cite{barbu1}),%
\begin{equation}
f_{\varepsilon }=\mathbf{1}_{\omega }\phi _{\varepsilon }e^{\frac{3}{2}%
s\alpha }.  \label{fcontrols}
\end{equation}%
Here, $(\phi _{\varepsilon },\theta _{\varepsilon })$ is the solution of the
adjoint system following:%
\begin{equation}
\begin{cases}
-\partial _{t}\phi _{\varepsilon }=\Delta \phi _{\varepsilon }+B\nabla \phi
_{\varepsilon }+\delta \theta _{\varepsilon } & \mathrm{in}\ Q, \\
-\partial _{t}\theta _{\varepsilon }=\Delta \theta _{\varepsilon }-\gamma
\theta _{\varepsilon }-\nabla \cdot \left( a\nabla \phi _{\varepsilon
}\right) & \mathrm{in}\ Q, \\
{\partial }_{\nu }\phi _{\varepsilon }=0,\partial _{\nu }\theta
_{\varepsilon }=0 & \mathrm{on}\ \Sigma , \\
(\phi _{\varepsilon },\theta _{\varepsilon })(x,T)=-\frac{1}{\varepsilon }%
(y_{\varepsilon },z_{\varepsilon })(x,T) & x\in \Omega ,%
\end{cases}
\label{adjointe}
\end{equation}%
where $(y_{\varepsilon },z_{\varepsilon })$ is the solution of %
\eqref{lineare} with $f=f_{\varepsilon }.$ By \eqref{lineare}, %
\eqref{fcontrols}, \eqref{adjointe}, and Proposition \ref{pro-obser}, it
follows that%
\begin{equation}
\iint_{Q_{\omega }}\left\vert \phi _{\varepsilon }\right\vert ^{2}e^{\frac{3%
}{2}s\alpha }dxdt+\frac{1}{\varepsilon }\left( \left\vert y(\cdot
,T)\right\vert _{2}^{2}+\left\vert z(\cdot ,T)\right\vert _{2}^{2}\right)
\leq e^{C\kappa }\left( \left\vert y_{0}\right\vert _{2}^{2}+\left\vert
z_{0}\right\vert _{2}^{2}\right) .  \label{null0}
\end{equation}%
We can simply get from \eqref{fcontrols} and \eqref{null0} that the control
function $f_{\varepsilon }$ satisfies%
\begin{equation*}
\left\Vert f_{\varepsilon }\right\Vert _{2}^{2}\leq e^{C\kappa }\left(
\left\vert y_{0}\right\vert _{2}^{2}+\left\vert z_{0}\right\vert
_{2}^{2}\right) .
\end{equation*}

Next we show that $f_{\varepsilon }$ can be taken in $L^{\infty }(Q)$. To
this end, let $\tau $ be a sufficiently small positive constant and let $%
\left\{ \tau _{j}\right\} _{j=0}^{M+1}$ be a finite increasing sequence such
that $0<\tau _{j}<\tau ,j=0,1,\ldots ,M,\tau _{M+1}=\tau .$ Let $\left\{
p_{i}\right\} _{i=0}^{M}$ be another finite increasing sequence such that $%
p_{0}=2,p_{M}>(N+2)/2$ and,
\begin{equation}
\left( \frac{N}{2}+1\right) \left( \frac{1}{p_{i}}-\frac{1}{p_{i+1}}\right) <%
\frac{1}{4}, \ i=0,1,\ldots, M-1.  \label{pipi}
\end{equation}
Set
\begin{equation*}
\alpha _{0}=\min_{\overline{\Omega }}\alpha =\frac{1-e^{2s\left\Vert \beta
\right\Vert _{C(\overline{\Omega })}}}{t(T-t)}.
\end{equation*}%
By \eqref{alpha},
\begin{equation*}
\alpha _{0}\leq \alpha \leq \frac{\alpha _{0}}{1+\eta (\lambda )}<0,
\end{equation*}%
where $\eta (\lambda )$ is defined by \eqref{eta}.

For each $i,i=0,1,\ldots ,M,M+1,$ define%
\begin{eqnarray*}
\zeta _{i}(x,t) &=&e^{(s+\tau _{i})\alpha _{0}}\phi _{\varepsilon }(x,T-t),
\\
\varrho _{i}(x,t) &=&e^{(s+\tau _{i})\alpha _{0}}\theta _{\varepsilon
}(x,T-t), \\
G_{i}(x,t) &=&\left[ \partial _{t}(e^{(s+\tau _{i})\alpha _{0}})\right] \phi
_{\varepsilon }(x,T-t), \\
H_{i}(x,t) &=&\left[ \partial _{t}(e^{(s+\tau _{i})\alpha _{0}})\right]
\theta _{\varepsilon }(x,T-t),
\end{eqnarray*}%
and
\begin{equation*}
\tilde{a}(x,t)=a(x,T-t),\tilde{B}(x,t)=B(x,T-t).
\end{equation*}%
Then for each $i,$ $(\zeta _{i},\varrho _{i})$ solves the following system:%
\begin{equation}
\begin{cases}
\partial _{t}\zeta _{i}-\Delta \zeta _{i}=\tilde{B}\nabla \zeta _{i}+\delta
\varrho _{i}+G_{i} & \mathrm{in}\ Q, \\
\partial _{t}\varrho _{i}-\Delta \varrho _{i}=-\gamma \varrho _{i}-\nabla
\cdot \left( \tilde{a}\nabla \zeta _{i}\right) +H_{i} & \mathrm{in}\ Q, \\
{\partial }_{\nu }\zeta _{i}=0,\partial _{\nu }\varrho _{i}=0 & \mathrm{on}\
\Sigma , \\
\zeta _{i}(x,0)=0,\varrho _{i}(x,0)=0 & x\in \Omega .%
\end{cases}
\label{adjointee}
\end{equation}%
Now we apply the $L^{p}$-$L^{q}$ estimate to the above system. By the
semigroup theory, the solution $(\zeta _{i},\varrho _{i})$, $i=1,2,\ldots
,M+1$, of \eqref{adjointee} can be represented as%
\begin{eqnarray}
\zeta _{i}(\cdot ,t) &=&\int_{0}^{t}e^{-(t-s)A}\left[ \tilde{B}\nabla \zeta
_{i}+\delta \varrho _{i}+G_{i}\right] (\cdot ,s)ds,  \label{zetai} \\
\varrho _{i}(\cdot ,t) &=&\int_{0}^{t}e^{-(t-s)A}\left[ -\gamma \varrho
_{i}-\nabla \cdot \left( \tilde{a}\nabla \zeta _{i}\right) +H_{i}\right]
(\cdot ,s)ds.  \label{rhoi}
\end{eqnarray}%
Firstly, by \eqref{semigroup1} to \eqref{zetai}, we have%
\begin{equation*}
\left\vert \zeta _{i}(\cdot ,t)\right\vert _{p_{i}}=C\int_{0}^{t}m(t-s)^{-%
\frac{N}{2}\left( \frac{1}{p_{i-1}}-\frac{1}{p_{i}}\right) }\left\vert
\left( \tilde{B}\nabla \zeta _{i}+\delta \varrho _{i}+G_{i}\right) (\cdot
,s)\right\vert _{p_{i-1}}ds,
\end{equation*}%
which can be estimated by Young's convolution inequality (see, e.g. \cite[p.3%
]{davies}) as
\begin{eqnarray}
\left\Vert \zeta _{i}\right\Vert _{p_{i}} &\leq &C\left( \left\Vert
B\right\Vert _{\infty }\left\Vert \nabla \zeta _{i}\right\Vert
_{p_{i-1}}+\left\Vert \varrho _{i}\right\Vert _{p_{i-1}}+\left\Vert
G_{i}\right\Vert _{p_{i-1}}\right)  \notag \\
&&\times \left[ \int_{0}^{T}m(t)^{-\frac{N}{2}\left( \frac{1}{p_{i-1}}-\frac{%
1}{p_{i}}\right) r_{i}}\right] ^{\frac{1}{r_{i}}}  \label{zzeta}
\end{eqnarray}%
where $r_{i}=1/[1-(1/p_{i-1})+(1/p_{i})]$. Similarly, applying %
\eqref{semigroup1} and \eqref{semigroup4} with $\varepsilon =\frac{1}{4}$ to %
\eqref{rhoi}, we have
\begin{eqnarray*}
\left\vert \varrho _{i}(\cdot ,t)\right\vert _{p_{i}}
&=&C\int_{0}^{t}m(t-s)^{-\frac{N}{2}\left( \frac{1}{p_{i-1}}-\frac{1}{p_{i}}%
\right) }\left\vert \left( -\gamma \varrho _{i}+H_{i}\right) (\cdot
,s)\right\vert _{p_{i-1}}ds \\
&&+C\left\Vert a\right\Vert _{\infty }\int_{0}^{t}m(t-s)^{-\frac{N}{2}\left(
\frac{1}{p_{i-1}}-\frac{1}{p_{i}}\right) -\frac{1}{2}-\frac{1}{4}}\left\vert
\nabla \zeta _{i}(\cdot ,s)\right\vert _{p_{i-1}}ds,
\end{eqnarray*}%
which can also be estimated by Young's convolution inequality as%
\begin{eqnarray}
\left\Vert \varrho _{i}\right\Vert _{p_{i}} &\leq &C\left( \left\Vert
\varrho _{i}\right\Vert _{p_{i-1}}+\left\Vert H_{i}\right\Vert
_{p_{i-1}}\right) \left[ \int_{0}^{T}m(t)^{-\frac{N}{2}\left( \frac{1}{%
p_{i-1}}-\frac{1}{p_{i}}\right) r_{i}}dt\right] ^{\frac{1}{r_{i}}}  \notag \\
&&+C\left\Vert a\right\Vert _{\infty }\left\Vert \nabla \zeta
_{i}\right\Vert _{p_{i-1}}\left[ \int_{0}^{T}m(t)^{\left( -\frac{N}{2}\left(
\frac{1}{p_{i-1}}-\frac{1}{p_{i}}\right) -\frac{3}{4}\right) r_{i}}dt\right]
^{\frac{1}{r_{i}}}.  \label{rrho}
\end{eqnarray}%
Owing to \eqref{pipi}, we also have
\begin{equation}
\left[ \int_{0}^{T}m(t)^{-\frac{N}{2}\left( \frac{1}{p_{i-1}}-\frac{1}{p_{i}}%
\right) r_{i}}dt\right] ^{\frac{1}{r_{i}}}\leq C\left[ (T+1)^{\frac{1}{r_{i}}%
}+T^{-\left( \frac{N}{2}+1\right) \left( \frac{1}{p_{i}}-\frac{1}{p_{i+1}}%
\right) +1}\right]  \label{mt}
\end{equation}%
and
\begin{equation}
\left[ \int_{0}^{T}m(t)^{\left( -\frac{N}{2}\left( \frac{1}{p_{i-1}}-\frac{1%
}{p_{i}}\right) -\frac{3}{4}\right) r_{i}}dt\right] ^{\frac{1}{r_{i}}}\leq C%
\left[ (T+1)^{\frac{1}{r_{i}}}+T^{-\left( \frac{N}{2}+1\right) \left( \frac{1%
}{p_{i}}-\frac{1}{p_{i+1}}\right) +\frac{1}{4}}\right] .  \label{ttt}
\end{equation}%
Secondly, we estimate the energy of solution $(\zeta _{i},\varrho _{i})$ to
get the following $L^{p_{i-1}}$-estimate%
\begin{equation}
\left\Vert \zeta _{i}\right\Vert _{p_{i-1}}+\left\Vert \varrho
_{i}\right\Vert _{p_{i-1}}+\left\Vert \nabla \zeta _{i}\right\Vert
_{p_{i-1}}\leq e^{C\kappa }\left( \left\Vert G_{i}\right\Vert
_{p_{i-1}}+\left\Vert H_{i}\right\Vert _{p_{i-1}}\right) .  \label{energg}
\end{equation}%
This inequality together with \eqref{zzeta}, \eqref{rrho}, \eqref{mt}, and %
\eqref{ttt} gives
\begin{equation}
\left\Vert \zeta _{i}\right\Vert _{p_{i}}+\left\Vert \varrho _{i}\right\Vert
_{p_{i}}\leq e^{C\kappa }\left( \left\Vert G_{i}\right\Vert
_{p_{i-1}}+\left\Vert H_{i}\right\Vert _{p_{i-1}}\right) .  \label{rozeta}
\end{equation}%
Since
\begin{equation}
\left\Vert G_{i}\right\Vert _{p_{i-1}}\leq CT\left\Vert \zeta
_{i-1}\right\Vert _{p_{i-1}}\text{ \ and }\left\Vert H_{i}\right\Vert
_{p_{i-1}}\leq CT\left\Vert \varrho _{i-1}\right\Vert _{p_{i-1}},
\label{ghi}
\end{equation}%
it follows from \eqref{rozeta} that
\begin{equation}
\left\Vert \zeta _{i}\right\Vert _{p_{i}}+\left\Vert \varrho _{i}\right\Vert
_{p_{i}}\leq e^{C_{i}\kappa }\left( \left\Vert \zeta _{i-1}\right\Vert
_{p_{i-1}}+\left\Vert \varrho _{i-1}\right\Vert _{p_{i-1}}\right) ,
\label{iteration}
\end{equation}%
where $C_{i}=C_{i}(\Omega ,\omega ),i=0,1,\ldots ,M,$ are positive
constants. The iteration inequality \eqref{iteration} from $0$ to $M$
implies that
\begin{equation}
\left\Vert \zeta _{M}\right\Vert _{p_{M}}+\left\Vert \varrho _{M}\right\Vert
_{p_{M}}\leq e^{C\kappa }\left( \left\Vert \zeta _{0}\right\Vert
_{2}+\left\Vert \varrho _{0}\right\Vert _{2}\right) .  \label{zetam}
\end{equation}%
By the definition of $\zeta _{0}$ and $\varrho _{0}$, we obtain from%
\eqref{null0} and \eqref{zetam} that
\begin{equation}
\left\Vert \zeta _{M}\right\Vert _{p_{M}}+\left\Vert \varrho _{M}\right\Vert
_{p_{M}}\leq e^{C\kappa }\left( \left\Vert y_{0}\right\Vert _{2}+\left\Vert
z_{0}\right\Vert _{2}\right) .  \label{zetam1}
\end{equation}%
Finally, we apply $L^{p_{M}}$-maximal regularity for the first equation of %
\eqref{adjointee} for $\zeta _{M+1}$ to get
\begin{eqnarray*}
&&\left\Vert \partial _{t}\zeta _{M+1}\right\Vert _{p_{M}}+\left\Vert \Delta
\zeta _{M+1}\right\Vert _{p_{M}}+\left\Vert \zeta _{M+1}\right\Vert _{p_{M}}
\\
&\leq &C\left( \left\Vert B\right\Vert _{\infty }\left\Vert \nabla \zeta
_{M+1}\right\Vert _{p_{M}}+\left\Vert \varrho _{M+1}\right\Vert
_{p_{M}}+\left\Vert G_{M+1}\right\Vert _{p_{M}}\right) .
\end{eqnarray*}%
This, by taking into account of \eqref{ghi} and \eqref{energg}, leads to%
\begin{equation*}
\left\Vert \zeta _{M+1}\right\Vert _{W_{p_{M}}^{2,1}(Q)}\leq e^{C\kappa
}\left( \left\Vert \zeta _{M}\right\Vert _{p_{M}}+\left\Vert \varrho
_{M}\right\Vert _{p_{M}}\right) .
\end{equation*}%
Hence, by the imbedding inequality (\cite[Lemma 3.3,Ch.II]{lady})
\begin{equation*}
\left\Vert \zeta _{M+1}\right\Vert _{C(\overline{Q})}\leq e^{C(1+T+\frac{1}{T%
})}\left\Vert \zeta _{M+1}\right\Vert _{W_{p_{M}}^{2,1}(Q)},
\end{equation*}%
for $p_{M}>(N+2)/2,$ and by \eqref{zetam1}, we get
\begin{equation*}
\left\Vert \zeta _{M+1}\right\Vert _{\infty }\leq e^{C\kappa }\left(
\left\Vert y_{0}\right\Vert _{2}+\left\Vert z_{0}\right\Vert _{2}\right) .
\end{equation*}%
That is
\begin{equation*}
\left\Vert \phi _{\varepsilon }e^{(s+\tau )\alpha _{0}}\right\Vert _{\infty
}\leq e^{C\kappa }\left( \left\Vert y_{0}\right\Vert _{2}+\left\Vert
z_{0}\right\Vert _{2}\right) .
\end{equation*}%
which together with \eqref{fcontrols} yields
\begin{equation*}
\left\Vert e^{\left[ -s(\frac{1}{2}-\eta (\lambda ))+\tau (1+\eta (\lambda ))%
\right] \alpha }f_{\varepsilon }\right\Vert _{\infty }\leq e^{C\kappa
}\left( \left\Vert y_{0}\right\Vert _{2}+\left\Vert z_{0}\right\Vert
_{2}\right) ,
\end{equation*}%
where $\eta (\lambda )$ is given by \eqref{eta}. This gives, by choosing $%
\tau $ small enough such that
\begin{equation*}
-s\left( \frac{1}{2}-\eta (\lambda )\right) +\tau (1+\eta (\lambda ))<0,
\end{equation*}%
that
\begin{equation*}
\left\Vert f_{\varepsilon }\right\Vert _{\infty }\leq e^{C\kappa }\left(
\left\Vert y_{0}\right\Vert _{2}+\left\Vert z_{0}\right\Vert _{2}\right) .
\end{equation*}%
The above inequality enables us to extract a subsequences of $f_{\varepsilon
}$, still denoted by itself, such that $f_{\varepsilon }\rightarrow f$
weakly in $L^{2}(Q)$, weakly$^{\ast }$ in $L^{\infty }(Q)$ as $\varepsilon
\rightarrow 0$. Let $(y_{\varepsilon },z_{\varepsilon })$ be the solution to
the system associated to $f_{\varepsilon }$. Then, by Proposition \ref{pro1}%
, we see that $y_{\varepsilon }$ and $z_{\varepsilon }$ are both bounded in $%
V^{1}(Q)$. Thus, there exist subsequences $y_{\varepsilon }$ and $%
z_{\varepsilon }$, still denoted by themselves, such that
\begin{equation*}
y_{\varepsilon }\rightarrow y,z_{\varepsilon }\rightarrow z\text{ \ weakly
in }V^{1}(Q)\text{; strongly in }L^{2}(Q)
\end{equation*}%
for $(y,z)\in V^{1}(Q)\cap C([0,T];L^{2}(\Omega ))$, which is the weak
solution of the system corresponding to $f\in L^{\infty }(Q)$, and $y(x,T)=0$
and $z(x,T)=0$ almost everywhere in $\Omega $. This completes the proof.
\hfill $\Box $

\section{Proof of Theorem \protect\ref{th-non}}

Let $(\overline{u},\overline{v})$ be a trajectory of the system \eqref{ks}
with the initial value $(\overline{u}_{0},\overline{v}_{0})$, which
satisfies \eqref{regularity}. Set $u=\overline{u}+y$, $v=\overline{v}+z$, $%
y_{0}=u_{0}-\overline{u}_{0}$, $z_{0}=v_{0}-\overline{v}_{0}$. Then, $(y,z)$
solves the following parabolic system
\begin{equation}
\begin{cases}
\partial _{t}y=\Delta y-\chi \nabla \cdot \left( y\nabla \overline{v}\right)
-\chi \nabla \cdot \left( (\overline{u}+y)\nabla z\right) +{\mathbf{1}}%
_{\omega }f & \mathrm{in}\ Q, \\
\partial _{t}z=\Delta z-\gamma z+\delta y & \mathrm{in}\ Q, \\
{\partial }_{\nu }y=0,\partial _{\nu }z=0 & \mathrm{on}\ \Sigma , \\
y(x,0)=y_{0}(x)\ \ z(x,0)=z_{0}(x) & x\in \Omega.%
\end{cases}
\label{eee}
\end{equation}%
The local exact controllability of the system \eqref{e} is equivalent to the
local null controllability of the system \eqref{eee}.

Let $K=\left\{ \eta \in L^{\infty }(Q)|\left\Vert \eta \right\Vert _{\infty
}\leq 1\right\} $. For each $\eta \in K$, we consider the following
linearized system%
\begin{equation}
\begin{cases}
\partial _{t}y=\Delta y-\nabla \cdot \left( By\right) -\nabla \cdot \left(
a_{\eta }\nabla z\right) +{\mathbf{1}}_{\omega }f & \mathrm{in}\ Q, \\
\partial _{t}z=\Delta z-\gamma z+\delta y & \mathrm{in}\ Q, \\
{\partial }_{\nu }y=0,\partial _{\nu }z=0 & \mathrm{on}\ \Sigma , \\
y(x,0)=y_{0}(x)\ \ z(x,0)=z_{0}(x) & x\in \Omega ,%
\end{cases}
\label{elinearized}
\end{equation}%
where $a_{\eta }=\chi (\overline{u}+\eta )$ and $B=\chi \nabla \overline{v}$%
. By \eqref{regularity}, we see that
\begin{equation*}
\text{ }a_{\eta }\in L^{\infty }(Q),\text{ \ }B\in L^{\infty }(Q)^{N}\text{
with }B\cdot \nu =0\text{ on }\Sigma .
\end{equation*}%
so system \eqref{elinearized} is casted into the exact framework of system %
\eqref{lineare}. Thus, we can apply Theorem \ref{th-lin} to obtain that for
each $\eta \in K$, there exists a pair $\left( (y,z),f\right) $which solves
system \eqref{elinearized} with $y(x,T)=0,z(x,T)=0$ almost everywhere in $%
\Omega $. Here and in what follows, we denote by $(y,z)$ the solution to
system \eqref{elinearized} corresponding to $f$ and $\eta $ if there is no
ambiguity. By \eqref{fcontrol}, we see that the control functions are
bounded as follows:%
\begin{equation}
\left\Vert f\right\Vert _{\infty }\leq e^{C\kappa _{0}}\left( \left\vert
y_{0}\right\vert _{2}+\left\vert z_{0}\right\vert _{2}\right) ,
\label{fcontrolin0}
\end{equation}%
where%
\begin{equation}
\kappa _{0}=c_{0}\left( 1+T+\frac{1}{T}\right) .  \label{k0}
\end{equation}%
By \eqref{thlinfty} of Proposition \ref{pro1} and \eqref{fcontrolin0}, we
have the following estimate%
\begin{equation}
\left\Vert y\right\Vert _{V^{1}(Q)}+\left\Vert z\right\Vert
_{V^{2}(Q)}+\left\Vert y\right\Vert _{\infty }+\left\Vert z\right\Vert
_{\infty }\leq e^{C\kappa _{0}}\left( \left\vert y_{0}\right\vert _{\infty
}+\left\Vert z_{0}\right\Vert _{W^{1,q_{N}}(\Omega )}\right) .
\label{rproof0}
\end{equation}%
For $\eta \in K$, define a multi-valued mapping $\Lambda :K\rightarrow
2^{L^{2}(Q)}$ by
\begin{equation*}
\Lambda (\eta )=\left\{
\begin{tabular}{ll}
$y\in L^{2}(Q)$ & $\left\vert
\begin{array}{c}
\exists f\text{ satisfying \eqref{fcontrolin0} such that }(y,z)\text{ is }
\\
\text{the solution to \eqref{elinearized} corresponding to }\eta \text{ and }f%
\text{, } \\
\text{and }y(x,T)=z(x,T)=0\text{ a.e. in }\Omega%
\end{array}%
\right. $%
\end{tabular}%
\right\} .
\end{equation*}%
We apply Kakutani's fixed-point theorem (\cite[p.7]{barbu1}) to the map $%
\Lambda $ to prove Theorem \ref{th-non}. First, it is clear that $K$ is a
convex subset of $L^{2}(Q)$. By the argument above, we see that $\Lambda
(\eta )$ is nonempty and convex for each $\eta \in K$. Moreover, by %
\eqref{rproof0}, $\Lambda (\eta )$ is bounded in $V^{1}(Q)$ for each $\eta
\in K$ and hence $\Lambda (\eta )$ is a compact subset of $L^{2}(Q)$ by the
Aubin-Lions lemma (\cite[p.17]{barbu1}).

Next, we show that $\Lambda $ is upper semi-continuous. To this purpose, let $%
\left\{ \eta _{n}\right\} _{n=1}^{\infty }$ be a sequence of functions in $K$
such that $\eta _{n}\rightarrow \eta $ strongly in $L^{2}(Q)$, and let $%
y_{n}\in \Lambda (\eta _{n})$ for each $n$. Then, by the definition of $%
\Lambda (\eta _{n})$, there exists $f_{n}$ for each $n$ such that $%
(y_{n},z_{n})$ solves the following system%
\begin{equation}
\begin{cases}
\partial _{t}y_{n}=\Delta y_{n}-\nabla \cdot \left( By_{n}\right) -\nabla
\cdot \left( a_{\eta _{n}}\nabla z_{n}\right) +{\mathbf{1}}_{\omega }f_{n} &
\mathrm{in}\ Q, \\
\partial _{t}z_{n}=\Delta z_{n}-\gamma z_{n}+\delta y_{n} & \mathrm{in}\ Q,
\\
{\partial }_{\nu }y_{n}=0,\partial _{\nu }z_{n}=0 & \mathrm{on}\ \Sigma , \\
y_{n}(x,0)=y_{0}(x)\ \ z_{n}(x,0)=z_{0}(x) & x\in \Omega ,%
\end{cases}
\label{eln}
\end{equation}%
and $y_{n}(x,T)=z_{n}(x,T)=0$ for $x\in \Omega $ almost everywhere.
Moreover, the control $f_{n}$ satisfies%
\begin{equation}
\left\Vert f_{n}\right\Vert _{\infty }\leq e^{C\kappa _{0}}\left( \left\vert
y_{0}\right\vert _{2}+\left\vert z_{0}\right\vert _{2}\right) .  \label{fcn}
\end{equation}%
By \eqref{fcn} and Proposition \ref{pro1}, we obtain
\begin{equation}
\left\Vert y_{n}\right\Vert _{V^{1}(Q)}+\left\Vert z_{n}\right\Vert
_{V^{2}(Q)}\leq e^{C\kappa }\left( \left\vert y_{0}\right\vert
_{2}+\left\Vert z_{0}\right\Vert _{W^{1,2}(\Omega )}\right) .  \label{ynzn}
\end{equation}%
By \eqref{fcn}, \eqref{ynzn} and applying the Aubin-Lions lemma again, we
can get $f\in L^{\infty }(Q)$, $y\in V^{1}(Q)$, $z\in V^{2}(Q)$ and the
subsequences of $f_{n}$, $y_{n}$, $z_{n}$, still denoted by themselves, such
that%
\begin{eqnarray*}
f_{n} &\rightarrow &f\text{ weakly}^* \hbox{ in }L^{\infty
}(Q),\text{ and weakly
in }L^{2}(Q)\text{;} \\
y_{n} &\rightarrow &y\text{ weakly in }V^{1}(Q)\text{, and strongly in }%
L^{2}(Q)\text{;} \\
z_{n} &\rightarrow &z\text{ weakly in }V^{2}(Q)\text{, and strongly in }%
L^{2}(0,T;H^{1}(\Omega )).
\end{eqnarray*}%
Passing to the limit as $n\rightarrow \infty $ in \eqref{eln}, we get that $%
(y,z)$ is a weak solution of \eqref{eln} corresponding to $\eta $.
We
claim that that $y\in \Lambda (\eta )$. Actually, let $Y_{n}=y_{n}-y$, $%
Z_{n}=z_{n}-z$, and $F_{n}=\mathbf{1}_{\omega }(f_{n}-f)$. Then $%
(Y_{n},Z_{n})$ solves the following system%
\begin{equation}
\begin{cases}
\partial _{t}Y_{n}=\Delta Y_{n}-\nabla \cdot \left( BY_{n}\right)  &  \\
\text{ \ \ \ \ \ \ }-\nabla \cdot \left[ a_{\eta _{n}}\nabla Z_{n}+(a_{\eta
_{n}}-a_{\eta })\nabla z\right] +F_{n} & \mathrm{in}\ Q, \\
\partial _{t}Z_{n}=\Delta Z_{n}-\gamma Z_{n}+\delta Y_{n} & \mathrm{in}\ Q,
\\
{\partial }_{\nu }Y_{n}=0,\partial _{\nu }Z_{n}=0 & \mathrm{on}\ \Sigma , \\
Y_{n}(x,0)=0\ \ Z_{n}(x,0)=0 & x\in \Omega .%
\end{cases}
\label{ynzncap}
\end{equation}%
Multiply the first equation of (\ref{ynzncap}) by $Y_{n}$, and integrate
over $\Omega $, to give
\begin{eqnarray}
\frac{d}{dt}\left\vert Y_{n}\right\vert _{2}^{2}+\left\vert \nabla
Y_{n}\right\vert _{2}^{2} &\leq &C\left\Vert B\right\Vert _{\infty
}^{2}\left\vert Y_{n}\right\vert _{2}^{2}+C\left\Vert a_{\eta
_{n}}\right\Vert _{\infty }^{2}\left\vert \nabla Z_{n}\right\vert _{2}^{2}
\notag \\
&&+C\int_{\Omega }\left\vert \eta _{n}-\eta \right\vert ^{2}\left\vert
\nabla z\right\vert ^{2}dx+C\int_{\Omega }F_{n}Y_{n}dx.  \label{ddd1}
\end{eqnarray}%
In the same way to the second equation of (\ref{ynzncap}), we have
\begin{equation}
\frac{d}{dt}\left\vert Z_{n}\right\vert _{2}^{2}+\left\vert \nabla
Z_{n}\right\vert _{2}^{2}+\gamma \left\vert Z_{n}\right\vert _{2}^{2}\leq
C\left\vert Y_{n}\right\vert _{2}^{2}.  \label{ddd2}
\end{equation}%
Differentiate $\left\vert \nabla Z_{n}\right\vert _{2}^{2}$ with respect to $%
t$ to get, from the second equation of (\ref{ynzncap}), that
\begin{equation}
\frac{d}{dt}\left\vert \nabla Z_{n}\right\vert _{2}^{2}+\left\vert \Delta
Z_{n}\right\vert _{2}^{2}+\gamma \left\vert \nabla Z_{n}\right\vert
_{2}^{2}\leq C\left\vert Y_{n}\right\vert _{2}^{2}.  \label{ddd3}
\end{equation}%
Since $\left\Vert a_{\eta _{n}}\right\Vert _{\infty }\leq C$, it follows
from \eqref{ddd1}-\eqref{ddd3} and Gronwall's lemma that
\begin{eqnarray}
&&\left\vert Y_{n}(\cdot ,t)\right\vert _{2}^{2}+\left\vert Z_{n}(\cdot
,t)\right\vert _{2}^{2}+\left\vert \nabla Z_{n}(\cdot ,t)\right\vert _{2}^{2}
\label{ynzngr} \\
&\leq &e^{C\left( 1+\left\Vert B\right\Vert _{\infty }^{2}\right) T}\left(
\int_{\Omega }\left\vert \eta _{n}-\eta \right\vert ^{2}\left\vert \nabla
z\right\vert ^{2}dx+\int_{\Omega }F_{n}Y_{n}dx\right) .  \notag
\end{eqnarray}%
On the other hand, since $(y,z)$ solves \eqref{elinearized}, by (ii) of
Proposition \ref{pro1}, we get that
\begin{equation*}
\left\Vert z\right\Vert _{W_{p}^{2,1}(Q)}\leq C\left( \left\vert
y_{0}\right\vert _{p}+\left\Vert z_{0}\right\Vert _{W^{2(1-\frac{1}{p}%
),p}(\Omega )}+\left\Vert \mathbf{1}_{\omega }f\right\Vert _{p}\right) ,
\end{equation*}%
which together with \eqref{fcontrolin0} implies
\begin{equation}
\left\Vert z\right\Vert _{W_{p}^{2,1}(Q)}\leq C\left( \left\vert
y_{0}\right\vert _{p}+\left\Vert z_{0}\right\Vert _{W^{2(1-\frac{1}{p}%
),p}(\Omega )}\right) .  \label{zzzz}
\end{equation}%
Since $W_{p}^{2,1}(Q)\hookrightarrow C^{1}(\overline{Q})$ for $p>N+2$ (\cite[%
Lemma 3.3, Ch II]{lady}), it follows from \eqref{zzzz} that
\begin{equation}
\left\Vert \nabla z\right\Vert _{C(\overline{Q})^{N}}\leq C\left( \left\vert
y_{0}\right\vert _{p}+\left\Vert z_{0}\right\Vert _{W^{2(1-\frac{1}{p}%
),p}(\Omega )}\right) .  \label{zzz0}
\end{equation}%
Since $\eta _{n}\rightarrow \eta $ strongly in $L^{2}(Q)$, $Y_{n}\rightarrow
0$ strongly in $L^{2}(Q)$, and $F_{n}\rightarrow 0$ weakly in $L^{2}(Q)$,
thus, by \eqref{zzz0}, we see that the right hand side of \eqref{ynzngr}
tends to $0$ as $n\rightarrow \infty $. Hence, $\left\vert Y_{n}(\cdot
,t)\right\vert _{2}\rightarrow 0$, $\left\vert Z_{n}(\cdot ,t)\right\vert
_{2}\rightarrow 0$ for all $t\in \lbrack 0,T]$. Since $%
y_{n}(x,T)=z_{n}(x,T)=0$ in $\Omega $ almost everywhere, we get that $%
y(x,T)=z(x,T)=0$ in $\Omega $ almost everywhere, which implies that $y\in
\Lambda (\eta )$. This shows that $\Lambda $ is upper semi-continuous.

Now it remains to show that $\Lambda (K)\subset K$. By Proposition \ref{pro1}%
, for any $y\in \Lambda (K)$,
\begin{equation*}
\left\Vert y\right\Vert _{\infty }\leq e^{c_{1}\kappa _{0}}\left( \left\vert
y_{0}\right\vert _{\infty }+\left\Vert z_{0}\right\Vert _{W^{2(1-\frac{1}{p}%
),p}(\Omega )}\right) ,
\end{equation*}%
where $c_{1}$ is a positive constant. Take $\delta =e^{-c_{1}\kappa
_{0}}$ such that if $\left\vert y_{0}\right\vert _{\infty
}+\left\Vert z_{0}\right\Vert _{W^{2(1-\frac{1}{p}),p}(\Omega )}\leq
\delta $ which is exactly \eqref{intialdata}, then $\left\Vert
y\right\Vert _{\infty }\leq 1$ and hence $\Lambda (K)\subset K$.
Therefore, the conditions of Kakutani's fixed point are satisfied,
that is, if the initial data $\left( u_{0},v_{0}\right) $ satisfies
\eqref{intialdata}, then there exists at
least one fixed point $y$, which together with $z$, is the solution of %
\eqref{e} corresponding with some control $f$ and satisfies $y(x,T)=0$ and $%
z(x,T)=0$ for $x\in \Omega $ almost everywhere. This completes the
proof.\hfill $\Box $


\bigskip

\begin{thebibliography}{99}
\bibitem{ammar} Ammar-Khodja, F., Benabdallah, A., Dupaix, C., Kostin, I.,
\textit{Controllability to the trajectories of phase-field models by one
control force}, SIAM J. Control Optim. 42 (2003), no. 5, 1661-1680.

\bibitem{ammar2} Ammar-Khodja, F., Benabdallah, A., Dupaix, C., \textit{%
Null-controllability of some reaction-diffusion systems with one control
force}, J. Math. Anal. Appl., 320 (2006), no. 2, 928-943.

\bibitem{ammar3} Ammar-Khodja, F., Benabdallah, A., Dupaix, C., Gonz\'{a}%
lez-Burgos, M., \textit{A generalization of the Kalman rank condition for
time-dependent coupled linear parabolic systems}, Differ. Equ. Appl., 1
(2009), no. 3, 427-457.

\bibitem{ammar4} Ammar-Khodja, F., Benabdallah, A., Dupaix, C., Gonz\'{a}%
lez-Burgos, M., \textit{A Kalman rank condition for the localized
distributed controllability of a class of linear parabolic systems}, J.
Evol. Equ., 9 (2009), 267-291.

\bibitem{ammar0} Ammar-Khodja, F., Benabdallah, A., Gonz\'{a}lez-Burgos, M.,
de Teresa, L., \textit{Recent results on the controllability of linear
coupled parabolic problems: a survey}, Math. Control Relat. Fields, 1
(2011), no. 3, 267-306.

\bibitem{arendt} Arendt, W., \textit{Semigroups and evolution equations:
functional calculus, regularity and kernel estimates, in Handbook of
Differential Equations: Evolutionary Equations, Vol 1, Elsevier, 2004}, 1-85.

\bibitem{barbu} Barbu, V., \textit{Controllability of parabolic and
Navier-Stokes equations, }Sci. Math. Japon., 56 (2002), no. 1,
143-211.

\bibitem{barbu1} Barbu, V., \textit{Analysis and Control of Nonlinear
Infinite-Dimensional Systems, }Academic Press, Boston, 1993.

\bibitem{biler} Biler, P., \textit{Local and global solvability of some
parabolic systems modelling chemotaxis, }Adv. Math. Sci. Appl., 8 (1998),
no. 2, 715-743.

\bibitem{gonzalez1} Bodart, O., Gonz\'{a}lez-Burgos, M., P\'{e}rez-Gac\'{\i}%
a, R., \textit{Insensitizing controls for a heat equation with a nonlinear
term involving state and the gradient}, Nonlinear Anal. 57 (2004),\ no. 5-6,
687-711.

\bibitem{boy} Boy, A., \textit{Analysis for a system of coupled
reaction-diffusion parabolic equations arising in biology,} Comput. Math.
Appl. 32 (1996), no. 4, 15-21.

\bibitem{coron} Coron, J.-M., \textit{Control and Nonlinearity,} AMS,
Providence, RI, 2007.

\bibitem{davies} Davies, E. B., \textit{Heat Kernels and Spectral Theory},
Cambridge University Press, Cambridge, 1989.

\bibitem{teresa} de Teresa, L., \textit{Insensitizing controls for a
semilinear heat equation,} Comm. Partial Differential Equations, 25 (2000),
no. 1-2, 39-72.

\bibitem{fernandez2} Fern\'{a}ndez-Cara, E., Gonz\'{a}lez-Burgos, M., de
Teresa, L., \textit{Boundary controllability of parabolic coupled equations,}
J. Funct. Anal., 259 (2010), no. 7, 1720-1758.

\bibitem{fernandez1} Fern\'{a}ndez-Cara, E., Guerrero, S., \textit{Global
Carleman inequalities for parabolic systems and applications to
controllability,} SIAM J. Control Optim., 45 (2006), no. 4, 1399-1446.

\bibitem{fernandez3} Fern\'{a}ndez-Cara, E., Zuazua, E., \textit{Null and
approximate controllability for weakly blowing up semilinear heat equations,}
Ann. Inst. H. Poincar\'{e} Anal. Non Lin\'{e}aire, 17 (2000), no. 5, 583-616.

\bibitem{fernandez4} Fern\'{a}ndez-Cara, E., Zuazua, E., \textit{The cost of
approximate controllability for heat equations: the linear case,} Adv.
Differential Equations, 5 (2000), no. 4-6, 465-514.

\bibitem{fursikov} Fursikov, A., Imanuvilov O. Yu., \textit{Controllability
of Evolution Equations,} Lecture Notes Series 34, Seoul National University,
Research Institute of Mathematics, Global Analysis Research Center, Seoul,
1996.

\bibitem{gajewski} Gajewski, H., Zacharias, K., \textit{Global behavior of a
reaction-diffusion system modelling chemotaxis, }Math. Nathr. 195 (1998),
77-114.

\bibitem{GigaSohr} Giga, Y., Sohr H., \textit{Abstract L}$^{p}$ \textit{%
estimates for the Cauchy problem with applications to the Navier-Stokes
equations in exterior domains, }J. Funct. Anal., 102 (1991), no. 1, 72-94.

\bibitem{gonzalez} Gonz\'{a}lez-Burgos, M., P\'{e}rez-Garc\'{\i}a, R.,
\textit{Controllability results for some nonlinear coupled parabolic systems
by one control force,} Asymptot. Anal. 46 (2006), no. 2, 123--162.

\bibitem{imanuvilov} Imanuvilov, O. Yu., Yamamoto, M., \textit{Carleman
inequalities for parabolic equations in Sobolev spaces of negative order and
exact controllability for semilinear parabolic equations,} Publ. Res. Inst.
Math. Sci., 39 (2003), no. 2, 227-274.

\bibitem{henry} Henry, D., \textit{Geometric Theory of Semilinear Parabolic
Equations}, Springer, Berlin, Heidelberg, 1987.

\bibitem{herrero} Herrero, M. A., Vel\'{a}zquez, J. J. L., \textit{A blow-up
mechanism for a chemotaxis model}, Ann. Scuola Norm. Sup. Pisa CI. Sci. (4),
24 (1997), no. 4, 633-683.

\bibitem{hillen} Hillen, T., Painter, K. J., \textit{A user's guide to PDE
models for chemotaxis}, J. Math. Biol. 58 (2009), 183-217.

\bibitem{horstmann1} Horstmann, D., \textit{From 1970 until present: the
Keller-Segel model in chemotaxis and its consequences I, }Jahresber.
Deutsch. Math.-Verein. 105 (2003), no. 3, 103--165.

\bibitem{horstmann} Horstmann, D., Winkler, M., \textit{Boundedness vs.
blow-up in a chemotaxis system, }J. Differential Equations, 215 (2005), no.
1, 52-107.

\bibitem{keller} Keller, E. F., Segel, L. A., \textit{Initiation of slime
mold aggregation viewed as an instability,} J. Theor. Biol., 26 (1970),
399-415.

\bibitem{lady} Ladyzhenskaja, O. A., Solonnikov, V. A., Ural'ceva N. N.%
\textit{, Linear and Quasi-linear Equations of Parabolic Type,} AMS,
Providence, RI, 1968.

\bibitem{lamberton} Lamberton, D., \textit{Equations d'\'{e}volution lin\'{e}%
aires associ\'{e}es \`{a} les semi-groupes de contractions dans les espaces L%
}$^{p}$, J. Funct. Anal., 72 (1987), no. 2, 252-262.

\bibitem{osaki} Osaki, K. and Yagi, A., \textit{Finite dimensional
attractors for one-dimensional Keller-Segel equations,} Funkcial. Ekva., 44
(2001), no. 3, 441-469.

\bibitem{rothe} Rothe, F., \textit{Global Solutions of Reaction-Diffusion
Systems,} LNM 1072, Springer-Verlag, 1984.

\bibitem{yagi1} Ryu, S.-U., Yagi, A., \textit{Optimal control of
Keller-Segel equations,} J. Math. Anal. Appl., 256 (2001), no. 1, 45-66.

\bibitem{yagi2} Yagi, A., \textit{Norm behavior of solutions to a parabolic
system of chemotaxis,} Math. Japon.,  45 (1997), no. 2, 241-265.

\bibitem{zhang} Wang, G., Zhang, L., \textit{Exact local controllability of
a one-control reaction-diffusion system,} J. Optim. Theory Appl., 131
(2006), no. 3, 453-467.
\end{thebibliography}
\end{document}